\documentclass[12pt,reqno]{amsart}
\usepackage[utf8]{inputenc}
\usepackage{fancyhdr}
\usepackage{amssymb}
\usepackage{amsthm}
\usepackage{amsmath}
\usepackage{geometry}
\usepackage{mathtools}
\usepackage{amsxtra}
\usepackage{mathrsfs}
\usepackage[all]{xy}
\usepackage{colonequals}
\usepackage{enumitem}
\usepackage{breqn}
\usepackage{caption}
\usepackage{placeins}

\theoremstyle{plain}
\newtheorem{thm}{Theorem}[section]
\newtheorem{thmx}{Theorem}

\newtheorem{lem}[thm]{Lemma} 
\newtheorem{prop}[thm]{Proposition}
\newtheorem{cor}[thm]{Corollary}

\theoremstyle{definition}
\newtheorem*{ack*}{Acknowledgment}
\theoremstyle{remark}

\newtheorem{rmk}{Remark}[section]

\DeclareMathOperator{\lcm}{lcm}

\newcommand{\C}{\mathbb C}

\newcommand{\N}{\mathbb N}
\newcommand{\PP}{\mathbb P}
\newcommand{\Q}{\mathbb Q}
\newcommand{\R}{\mathbb R}
\newcommand{\Z}{\mathbb Z}

\newcommand{\bdd}{\mathbf d}
\newcommand{\bdm}{\mathbf m}
\newcommand{\bdn}{\mathbf n}

\newcommand{\Aa}{\mathcal{A}}
\newcommand{\Bb}{\mathcal{B}}
\newcommand{\Cc}{\mathcal{C}}

\newcommand{\Ee}{\mathcal{E}}

\newcommand{\Hh}{\mathcal{H}}

\newcommand{\Nn}{\mathcal{N}}
\newcommand{\Pp}{\mathcal{P}}

\newcommand{\Rr}{\mathcal{R}}
\newcommand{\Ss}{\mathcal{S}}
\newcommand{\Tt}{\mathcal{T}}

\newcommand{\Ww}{\mathcal{W}}

\newcommand{\psmod}[1]{\,(\textup{\text{mod}}\,{#1})}

\input cyracc.def

\usepackage[hidelinks]{hyperref}
\makeatletter
\@namedef{subjclassname@2020}{\textup{2020} Mathematics Subject Classification}
\makeatother

\begin{document}
\title{The Hardy--Ramanujan inequality for sifted sets and its applications}
\author{Kai (Steve) Fan}
\address{Department of Mathematics\\ University of Georgia\\ Athens, GA 30602}
\email{Steve.Fan@uga.edu}
\setlength{\footskip}{6mm}
\begin{abstract}
The well-known Hardy--Ramanujan inequality states that if $\omega(n)$ denotes the number of distinct prime factors of a positive integer $n$, then there is an absolute constant $C>0$ such that uniformly for $x\ge2$ and $k\in\N$,
\[\#\{n\le x\colon\omega(n)=k\}\ll\frac{x(\log\log x+C)^{k-1}}{(k-1)!\log x}.\]
A myriad of generalizations and variations of this inequality have been discovered. In this paper, we establish a weighted version of this inequality for sifted sets, which generalizes an earlier result of Hal\'asz and implies Timofeev's theorems on shifted primes. We then explore its applications to a variety of intriguing problems, such as large deviations of $\omega$ on subsets of integers, the Erd\H{o}s multiplication table problem, divisors of shifted primes, and the image of the Carmichael $\lambda$-function. Building on the same circle of ideas, we also generalize Troupe's result on the normal order of $\omega(s(n))$ for the sum-of-proper-divisors function $s(n)$, confirming for the first time the weighted version of a special case of a 1992 conjecture by Erd\H{o}s, Granville, Pomerance, and Spiro.
\end{abstract}

\maketitle
\section{Introduction}\label{S:Intro}
For every $n\in\N$, let $\omega(n)$ count the number of distinct prime factors of $n$. The fundamental distribution of $\omega(n)$ has been known for almost a century. One of the earliest results is an
inequality due to Hardy and Ramanujan \cite[Lemma B]{HR17}, which states that there exists an absolute constant $C>0$ such that
\begin{equation}\label{eq:HR17}
\#\{n\le x\colon\omega(n)=k\}\ll\frac{x(\log\log x+C)^{k-1}}{(k-1)!\log x}
\end{equation}
uniformly for $x\ge2$ and $k\in\N$.\footnote{One may actually take $C=\sup_{t\ge2}\left|\sum_{p\le t}1/p-\log\log t\right|$.} This inequality was developed by Hardy and Ramanujan in their investigation of the question ``What is the normal degree of compositeness of a number $n$?" When the ``degree of compositeness" is measured by $\omega(n)$, this question asks essentially about the normal order of $\omega(n)$. The uniformity of \eqref{eq:HR17} in $k$, together with the standard estimates for the Poisson distribution $\text{Pois}(\lambda )$ applied with $\lambda=\log\log x+O(1)$, enabled Hardy and Ramanujan to deduce that $\omega(n)$ has normal order $\log\log n$ \cite[Theorem B]{HR17}, in the sense that for a given $\epsilon>0$, we have $|\omega(n)-\log\log n|<\epsilon\log\log n$ for all but $o(x)$ of $n\le x$. This result may be viewed as a precursor of the celebrated Erd\H{o}s--Kac theorem \cite{EK40}, which asserts that the limiting distribution of $\omega(n)$ is Gaussian with mean $\log\log n$ and variance $\log\log n$.
More precisely, it states that for every $T\in\R$ we have
\[\lim_{x\to\infty}\frac{1}{x}\cdot\#\left\{n\le x\colon \frac{\omega(n)-\log\log n}{\sqrt{\log\log n}}\le T\right\}=\frac{1}{\sqrt{2\pi}}\int_{-\infty}^{T}e^{-t^2/2}\,dt.\]
In a similar vein, Hardy and Ramanujan also showed that $\Omega(n)$, the total number of prime factors of $n$ counted with multiplicity, has the same normal order \cite[Theorem C]{HR17}, but the analogue to \eqref{eq:HR17} for $\Omega(n)$ \cite[Lemma C]{HR17} turns out to be much more complicated.

There have been numerous generalizations and variations of \eqref{eq:HR17} since the original work of Hardy and Ramanujan. For instance, Hal\'asz \cite{Hala72} derived a version of the Hardy--Ramanujan inequality on the number of prime factors of $n$ belonging to a given set of primes. An analogue of \eqref{eq:HR17} on the number of prime factors of shifted primes $p-1$ was obtained by Erd\H{o}s \cite{Erd35} and later refined by Timofeev \cite{Tim95}, the former of which had foreshadowed an analogous Erd\H{o}s--Kac law for $\omega(p-1)$ later confirmed by Halberstam \cite{Halb56}. It is also natural to count integers with a prescribed number of prime factors in short intervals and in arithmetic progressions. Relevant results in this direction can be found in \cite{CCT11,Ora82,Spiro85,Tud00,WW77}. In a recent series of developments, Goudout \cite{Goud17} proved a Hardy--Ramanujan type inequality for the joint distribution of values of $\omega$ at two linear polynomial arguments. His result was later generalized by Tenenbaum \cite{Ten18} to an inequality for the joint distribution of values of $\omega$ at several irreducible polynomial arguments. Adapting the method of Hardy and Ramanujan, Pollack \cite{Poll20a} proved a weighted version of \eqref{eq:HR17} where nonnegative multiplicative functions serve as the weights. Prior to Pollack's work, Tenenbaum \cite{Ten00} had obtained a similar result in his study of the joint distribution of large prime factors of $n$. More recently, Ford \cite{Ford21,Ford22} considered analogues of \eqref{eq:HR17} for integers with a prescribed number of prime factors in disjoint sets and for integers in a subset subject to a certain sieve condition, the latter of which can particularly be applied to shifted primes. 

In the present paper, we establish a version of the Hardy--Ramanujan inequality for nonnegative multiplicative functions over a sifted set. In particular, it generalizes the main theorem in \cite{Hala72} and implies \cite[Theorems 1, 2]{Tim95} for $k$ in the comparable ranges.

Recall that an arithmetic function $f\colon\N\to\R$ is said to be {\it multiplicative} if $f(1)=1$ and $f(mn)=f(m)f(n)$ whenever $\gcd(m,n)=1$. We shall be concerned with the class $\mathscr{M}(A_1,A_2)$ of multiplicative functions $f\colon \N\to\R_{\ge0}$, where $A_1>0$ is an arbitrary constant, and $A_2\colon\R_{>0}\to\R_{>0}$ is any given function, satisfying the following two conditions:
\begin{enumerate}[label=(\roman*), leftmargin=*]
	\item $f(n) \le A_1^{\Omega(n)}$ for all $n\in\N$.
	\item Given every $\epsilon > 0$, one has $f(n)\le A_2(\epsilon)n^{\epsilon}$ for all $n\in\N$.
\end{enumerate}
Pollack \cite[Theorem 1.1]{Poll20b} proved the following variant of Shiu's theorem \cite[Theorem 1]{Shiu80} for multiplicative functions in the class $\mathscr{M}(A_1,A_2)$ over a sifted set.
\begin{thmx}\label{thmx:Poll20b}
Let $\alpha\in(0,1)$, $v\in\Z_{\ge0}$, $A_1>0$, $A_2\colon\R_{>0}\to\R_{>0}$, and $f\in\mathscr{M}(A_1,A_2)$. Let $2\le x^{\alpha}\le y\le x$, and let $\Nn\colonequals\N\cap(x-y,x]$ be the subset of positive integers $n$ such that $n\psmod{p}\notin\Ee_p\subseteq(\Z/p\Z)^{\times}$ for every prime $p$, where $\#\Ee_p=\nu(p)\le v$. Then
\[\sum_{n\in\Nn}f(n)\ll_{A_1,A_2,v,\alpha}\frac{y}{\log x}\exp\left(\sum_{p\le x}\frac{f(p)-\nu(p)}{p}\right).\]
\end{thmx}
Theorem \ref{thmx:Poll20b} may be thought of as Brun’s upper bound sieve inequality with a multiplicative weight attached. We will make heavy use of Theorem \ref{thmx:Poll20b} in the development of our version of the Hardy--Ramanujan inequality, which is in fact partly inspired by this result.

To state our results in a coherent manner, we need to introduce some more notation. Let $\PP$ be the set of primes and put $\pi(x)\colonequals\#(\PP\cap[1,x])$. In the sequel we shall reserve the letters $p$ and $q$ for primes. For any subset $E\subseteq\PP$, let $E^c\colonequals\PP\setminus E$, and define
\begin{align*}
\omega(n,E)&\colonequals\#\{p\in E\colon p\mid n\},\\
\Omega(n,E)&\colonequals\#\{(p,\ell)\in E\times\N\colon p^{\ell}\mid n\}.
\end{align*}
These functions have been extensively studied in the literature, including the aforementioned works \cite{Ford21,Hala72,Tim95}. In addition, we set 
\[M_f(x,E)\colonequals\sum_{\substack{p\le x\\p\in E}}\frac{f(p)}{p}\]
for any function $f\colon\PP\to\R_{\ge0}$, with the abbreviations that  $M(x,E)\colonequals M_f(x,E)$ when $f(p)=1$ identically and that $M_f(x)\colonequals M_f(x,\PP)$. Our version of the Hardy--Ramanujan inequality can be formulated as the following theorem.
\begin{thm}\label{thm:HRS}
Let $x\ge2$, $p_0\in\PP$, $v\in\Z_{\ge0}$, $\alpha_1>0$, and $\alpha_2\in(0,p_0)$. Let $f\in\mathscr{M}(A_1,A_2)$ with some $A_1>0$ and $A_2\colon\R_{>0}\to\R_{>0}$, and denote by $\Ss$ the set of positive integers $n\le x$ such that $n\psmod{p}\notin\Ee_p\subseteq(\Z/p\Z)^{\times}$ for every prime $p$, where $\#\Ee_p=\nu(p)\le v$. Let $E\subseteq\PP$ with $\min E\ge p_0$, and put $g\in\{\omega,\Omega\}$ and $\beta\colonequals \alpha_11_{g=\omega}+\alpha_21_{g=\Omega}$. Then 
\[\sum_{\substack{n\in\Ss\\g(n,E)=k}}f(n)\ll_{A_1,A_2,p_01_{g=\Omega},v,\beta}\frac{x(M_f(x,E)+O(1))^{k}}{k!\log x}e^{M_{f}(x,E^c)-M_{\nu}(x)}\]
for all $0\le k\le \beta M_f(x,E)$. Furthermore, we have in the case $E=\PP$ that
\[\sum_{\substack{n\in\Ss\\g(n)=k}}f(n)\ll_{A_1,A_2,p_01_{g=\Omega},v,\beta}\frac{x(M_f(x)+O(1))^{k-1}}{(k-1)!\log x}e^{-M_{\nu}(x)}\]
for all $1\le k\le \beta M_f(x)$. The terms $O(1)$ in both estimates depend at most on $A_1,A_2,v$.
\end{thm}

The case $k=0$ of Theorem \ref{thm:HRS} follows instantly from the case $y=x$ of Theorem \ref{thmx:Poll20b} applied to $f1_{(n,E)=1}$. Conversely, one recovers the latter from the former with $g=\omega$, $E=\PP\cap(x,\infty)$, and $k=0$. The proof of the case where $k\ge1$ and $M_f(x,E)\ge1/\beta$ will be given in Section \ref{S:HRS}.

The following corollary of Theorem \ref{thm:HRS} shows that under additional hypotheses, $g(n,E)$ has the expected normal order $M_f(x,E)$ with respect to the natural probability measure induced by $f$ on $\Ss$ defined, for every $m\in\Ss$, by
\[\text{Prob}(\bdm=m)=\frac{f(m)}{\sum_{n\in\Ss}f(n)}.\] 

\begin{cor}\label{cor:WNO}
Assume the hypotheses of Theorem \ref{thm:HRS}. Let $D\ge1$ and $\Pp$ the set of primes $p$ with $\nu(p)>0$. Suppose also that $\max\Pp\le D^{1/s_v}$, where $s_v\colonequals 1+2/(e^{0.53/v}-1)$, that $\nu(p)<p-1$ for all $p\in\Pp$, that $M_f(x,E)\ge c_0>0$ whenever $x\ge x_0\ge2$, that 
\begin{equation}\label{eq:sumf(p)}
\sum_{p\le t}f(p)\ge B\frac{t}{\log 2t}
\end{equation}
for all $t\in[x^{\theta},x]$, where $B>0$ and $\theta\in(1/2,1)$ are constants, and that 
\begin{equation}\label{eq:fBV}
\sum_{\substack{d\le D\\d\mid\Pp}}\nu(d)\max_{(a,d)=1}\left|\sum_{\substack{n\le x\\n\equiv a\psmod{d}}}f(n)-\frac{F(d)}{\varphi(d)}\sum_{n\le x}f(n)\right|=o\left(\frac{x}{\log x}e^{M_f(x)-M_{\nu}(x)}\right)
\end{equation}
as $x\to\infty$, where
\[F(d)\colonequals\prod_{p\mid d}\left(\sum_{\ell\ge0}\frac{f(p^{\ell})}{p^{\ell}}\right)^{-1}.\]
Then for any $x\ge x_0$ and $0<\lambda\le\sqrt{M_f(x,E)}/2$, we have
\[\left(\sum_{n\in\Ss}f(n)\right)^{-1}\sum_{\substack{n\in\Ss\\|g(n,E)-M_f(x,E)|\ge\lambda \sqrt{M_f(x,E)}}}f(n)\ll\lambda^{-1}\exp\left(-\frac{\lambda^2}{2}+O\left(\frac{\lambda^3}{\sqrt{M_f(x,E)}}\right)\right),\]
where the implied constant in ``$\ll$" may depend on $A_1,A_2,B,c_0,v,x_0,\theta$ as well as the implied constant in \eqref{eq:fBV}, while the implied constant in the big O is absolute.
\end{cor}

Corollary \ref{cor:WNO} quantifies the deviations of $g(n,E)$ on $\Ss$ from its weighted mean $M_f(x,E)$. The assumption \eqref{eq:fBV} indicates that 
\[\sum_{\substack{n\le x\\n\equiv a\psmod{d}}}f(n)\approx\frac{1}{\varphi(d)}\sum_{\substack{n\le x\\(n,d)=1}}f(n)\approx\frac{F(d)}{\varphi(d)}\sum_{n\le x}f(n)\]
on average up to the level $D$, where $F(d)$ may be thought of as the probability that a randomly chosen integer $\bdn\in\N$ subject to the natural probability measure induced by $f$ on $\N$ is relatively prime to $d$. Besides, the assumption that $\nu(p)<p-1$ for all $p\in\Pp$ rules out some extreme cases such as when $\nu(p)=p-1$ for some $p\in\Pp$ and $f(p^{\ell})=0$ for all $\ell\ge1$. It is also possible to show under stronger conditions that the limiting distribution of $g(n,E)$ is Gaussian with mean $M_f(x,E)$ and standard variance $\sqrt{M_f(x,E)}$, in the sense that 
\[\lim_{x\to\infty}\left(\sum_{n\in\Ss}f(n)\right)^{-1}\sum_{\substack{n\in\Ss\\g(n,E)\le M_f(x,E)+T\sqrt{M_f(x,E)}}}f(n)=\frac{1}{\sqrt{2\pi}}\int_{-\infty}^{T}e^{-t^2/2}\,dt\]
for every $T\in\R$. Such results are often referred to as ``weighted Erd\H{o}s--Kac theorems". The interested reader may consult \cite{EG24,Ell12,Ell15,Fan25a,KMS22,Ten17,Ten20} for the special case $\Ss=[1,x]$.

Our next result is a concrete example, showing that under certain conditions, the function $g(n,E)$, with $n$ running over the set of positive integers representable by a primitive, positive-definite binary quadratic form, has relative normal order $M(x,E)$. Recall that a binary quadratic form $F(X,Y)=aX^2+bXY+cY^2$ is said to be {\it primitive} if $\gcd(a,b,c)=1$. It is positive-definite if and only if $a,c>0$ and its discriminant $\Delta_F\colonequals b^2-4ac$ is negative. 
\begin{cor}\label{cor:QFNO}
Let $p_0\in\PP$, $\alpha_1>1$ and $\alpha_2\in(1,p_0)$, and set $\beta=\alpha_11_{g=\omega}+\alpha_21_{g=\Omega}$, $g\in\{\omega,\Omega\}$ and $\theta=1/(1+4\log\beta)$. Let $F(X,Y)=aX^2+bXY+cY^2\in\Z[X,Y]$ be a primitive, positive-definite binary quadratic form of discriminant $\Delta$, and put $d=\gcd(a,c)$ and $\Aa_F\colonequals F(\Z^2)\cap[1,x]$. If $E\subseteq\PP$ is a subset such that $\min E\ge p_0$ and $(\Delta/p)=1$ for $p\in E$, then
\begin{align*}
&\hspace*{5mm}\frac{1}{\#\Aa_F}\cdot\#\left\{n\in\Aa_F\colon|g(n,E)-M(x,E)|\ge\lambda\sqrt{M(x,E)}\right\}\\
&\ll_{d,p_0,\beta,\Delta}\lambda^{-1}\exp\left(-\frac{((1-\theta)\lambda)^2}{2}+O\left(\frac{\lambda^3}{\sqrt{M(x,E)}}\right)\right)
\end{align*}
for all $x\ge p_0$ and $0<\lambda\le\sqrt{M(x,E)}/2$, where the implied constant in the big O is absolute. In addition, if $F$ represents arbitrarily large numbers $n\not\equiv k\psmod{2}$, where $k\in\Z\setminus\{0\}$, then the same estimate holds with $\Bb_F\colonequals\Aa_F\cap(\PP-k)=\{n\in\Aa_F\colon n+k\in\PP\}$ in place of $\Aa_F$ with the implied constant in ``$\ll$" also dependent on $k$.
\end{cor}

Similarly, we have the following analogue for linear forms $ap+b$ in a prime variable $p$.
\begin{cor}\label{cor:SPNO}
Let $a\in\N$ and $b\in\Z\setminus\{0\}$ with $\gcd(a,b)=1$. Let $B>0$, $p_0\in\PP$, and $f\in\mathscr{M}(A_1,A_2)$ with some $A_1>0$ and $A_2\colon\R_{>0}\to\R_{>0}$. Suppose that $f(p)\ge B$ for all $p\ge p_0$ and that $f(2^{\ell})\ge B$ for some $\ell\in\N$ when $2\nmid ab$. If $E\subseteq\PP$ is a nonempty subset such that $M_f(x,E)\ge c_0>0$ whenever $x\ge x_0\ge\max\{2,(p_0-b)/a\}$, and if $g\in\{\omega,\Omega\}$, then 
\begin{align*}
&\hspace*{5mm}\left(\sum_{p\le x}f(ap+b)\right)^{-1}\sum_{\substack{p\le x\\|g(ap+b,E)-M_f(x,E)|\ge\lambda \sqrt{M_f(x,E)}}}f(ap+b)\\
&\ll_{A_1,A_2,B,a,b,c_0,\ell1_{2\nmid ab},x_0}\lambda^{-1}\exp\left(-\frac{\lambda^2}{2}+O\left(\frac{\lambda^3}{\sqrt{M_f(x,E)}}\right)\right)
\end{align*}
for all $x\ge x_0$ and $0<\lambda\le\sqrt{M(x,E)}/2$, where the implied constant in the big O is absolute.
\end{cor}

For each $N\in\N$, let $A(N)$ record the number of distinct entries in the 
$N\times N$ multiplication table. In other words, $A(N)=\#\{ab\colon a,b\le N\}$. Erd\H{o}s \cite{Erd60} showed that $A(N)=N^2/(\log N)^{\eta_0+o(1)}$ for sufficiently large $N$, where 
\begin{equation}\label{eq:eta_0}
\eta_0\colonequals 1-\frac{1+\log\log 2}{\log 2}=0.0860713...
\end{equation}
is the Erd\H{o}s--Tenenbaum--Ford constant. In his seminal work \cite{Ford08} on integers with a divisor in a given interval, Ford determined the correct order of magnitude for $A(N)$:
\begin{equation}\label{eq:A(N)}
A(N)\asymp\frac{N^2}{(\log N)^{\eta_0}(\log\log N)^{3/2}}.
\end{equation}
Koukoulopoulos \cite{Kou10} later obtained an analogue for the multiplication table of shifted primes $p-s$ for any fixed $s\in\Z\setminus\{0\}$, showing that
\begin{equation}\label{eq:P_s(N)}
A(N;P_s)\colonequals\#\{ab\in\PP-s\colon a,b\le N\}\asymp_s\frac{A(N)}{\log N}.
\end{equation}
By a short and crude argument, we obtain readily from Theorem \ref{thm:HRS} the following result on the distinct entries in a ``sifted" multiplication table.

\begin{cor}\label{cor:ErdosMT}
Let $f$ and $\Ss$ be as in Theorem \ref{thm:HRS}, and let $R_f(x)\colonequals (M_f(x)\log2)/\log\log x$ and $\Ss^{\ast}\colonequals\{n\in\Ss\colon n=ab\text{~with some~}a,b\le\sqrt{x}\}$. Suppose that there exist $c_0>0$ and $x_0\ge3$ such that $M_f(x)\ge c_0$ and $R_f(x)<1$ whenever $x\ge x_0$. Then for all $x\ge x_0$ we have that
\[\sum_{n\in\Ss^{\ast}}f(n)\ll_{A_1,A_2,c_0,v,x_0}\left(1+\frac{4^{M_f(x)}\sqrt{M_f(x)}}{\log x}\right)\frac{x}{(\log x)\sqrt{M_f(x)}}e^{2(1-\log 2)M_f(x)-M_{\nu}(x)}\]
when $R_f(x)\le1/2$, and that
\[\sum_{n\in\Ss^{\ast}}f(n)\ll_{A_1,A_2,c_0,v,x_0}\left(\frac{1}{1-R_f(x)}+\frac{1}{\sqrt{2R_f(x)-1}}\right)\frac{x}{(\log x)\sqrt{M_f(x)}}e^{(1-Q(1/R_f(x)))M_f(x)-M_{\nu}(x)}\]
when $1/2<R_f(x)<1$, where we define, for each $y\ge0$, 
\[Q(y)\colonequals\int_{1}^{y}\log t\,dt=y\log y-y+1.\]
\end{cor}

Despite its great generality, Corollary \ref{cor:ErdosMT} hardly captures the true order of magnitude (owing to the crude argument used in its proof). It is easy to see that for both $A(N)$ and $A(N;P_s)$, it yields the upper bounds which are only a factor of $\log\log N$ off their true orders of magnitude exhibited in \eqref{eq:A(N)} and \eqref{eq:P_s(N)}, respectively. More generally, for any $k$ distinct non-constant, irreducible polynomials $Q_1,...,Q_k\in\Z[X]$ with positive leading coefficients such that $Q_0=Q_1\cdots Q_k$ has no fixed prime factors, we have 
\[\sum_{n\in\Ss^{\ast}}1_{\PP}(Q_1(n))\cdots 1_{\PP}(Q_k(n))\ll_{Q_0}\frac{N^2}{(\log N)^{k+\eta_0}\sqrt{\log\log N}}\]
with $\Ss=\N\cap[1,N^2]$, which we expect to be off its true order of magnitude by a factor of $\log\log N$. On the other hand, taking $f(n)=r(n)/4$ instead, where $r(n)\colonequals\#\left\{(X,Y)\in\Z^2\colon X^2+Y^2=n\right\}$, we have analogously
\[\sum_{n\in\Ss^{\ast}}r(n)\ll\frac{N^2}{(\log N)^{\eta_0}\sqrt{\log\log N}},\]
whereas taking $f(n)=1_{r(n)>0}$, the characteristic function of sums of two squares, leads to
\[\#\{ab\colon a,b\le N\text{~and~}r(ab)>0\}\ll\frac{N^2}{(\log N)^{\log 2}\sqrt{\log\log N}}.\]
We hope to revisit these and related problems on a future occasion.

The proofs of Corollaries \ref{cor:WNO}--\ref{cor:ErdosMT} will be presented in Section \ref{S:NOg(n,E)}.

For each $n\in\N$, denote by $s(n)$ the sum of the proper divisors of $n$, namely,
\[s(n)\colonequals\sigma(n)-n=\sum_{d\mid n,\,d<n}d.\]
A 1992 conjecture of Erd\H{o}s, Granville, Pomerance, and Spiro \cite[Conjecture 4]{EGPS90} states that if $\Aa\subseteq\N$ has natural density zero, so does $s^{-1}(\mathcal{A})$. A few special cases of the EGPS conjecture have since been verified; see \cite{BCDDT24,Poll14, Poll15,PPT18,Tro15,Tro20}. The following is a generalization of a result of Troupe \cite[Theorem 1.3]{Tro15} where $\Aa$ consists of $n\in\N$ with abnormally many prime factors.
\begin{thm}\label{thm:EGPSTro}
Let $A_1>0$ and $A_2\colon\R_{>0}\to\R_{>0}$. Suppose that $f\in\mathscr{M}(A_1,A_2)$ satisfies \eqref{eq:sumf(p)} for sufficiently large $t$ with a constant $B>0$. For $\lambda=c_0\sqrt{\log_4 x}$ with any constant $c_0>2$, where $\log_4x\colonequals\log\log\log\log x$, we have
\[\lim_{x\to\infty}\left(\sum_{n\le x}f(n)\right)^{-1}\sum_{\substack{n\le x\\|\omega(s(n))-\log\log x|\ge\lambda\sqrt{\log\log x}}}f(n)=0.\]
\end{thm}

Theorem \ref{thm:EGPSTro} appears to be the first addressing weighted versions of the EGPS conjecture. While Troupe's argument is based on the second moment method of Tur\'{a}n, we will prove Theorem \ref{thm:EGPSTro} in Section \ref{S:EGPSTro} by combining the ideas from Section \ref{S:NOg(n,E)} with those from \cite{Poll14,Ten18}.

Along with the machinery developed in Section \ref{S:NOg(n,E)}, Theorem \ref{thm:HRS} also enables one to study shifted primes possessing a large shifted-prime divisor.
\begin{thm}\label{thm:LMPshifted}
For $a\in\Z\setminus\{0\}$, $u\in\N$ and $v\in\Z\setminus\{-au\}$, define
\[P_{a,u,v}(x,y)\colonequals\#\left\{p\le x\colon up+v\emph{~has a shifted-prime divisor~}q-a>y\right\}.\]
Then we have
\begin{equation}\label{eq:P_{a,u,v}(x,y)}
P_{a,u,v}(x,y)\ll_{a,u,v}\frac{\pi(x)}{(\log y)^{\eta_0}\sqrt{\log\log y}}
\end{equation}
for all $x,y\ge 3$, where $\eta_0$ is defined by \eqref{eq:eta_0}.
\end{thm}

Theorem \ref{thm:LMPshifted}, the proof of which can be found in Section \ref{S:SPlarge}, supplies an analogue to \cite[Theorem 1.2]{MPP17} on numbers with a large shifted-prime divisor by sharpening the bound $O(\pi(x)/(\log y)^{c})$ with some constant $c>0$ as displayed in \cite[Theorem 3]{LMP15}. On a related note, the distribution of the counting function of shifted-prime divisors on shifted primes was recently investigated by the author \cite{Fan25b}. In particular, an easy application of Brun's upper bound sieve as employed in the estimation of  \cite[Eq.\,(11)]{Fan25b} shows that under the hypothesis of Theorem \ref{thm:LMPshifted}, we have
\[P_{a,u,v}(x,y)\le\sum_{y<q-a\le ux+v}\sum_{\substack{p\le x\\(q-a)\mid (up+v)}}1\ll_{a,u,v}\frac{\pi(x)\max\{\log(x/y),1\}}{\log x}\]
for all $x,y\ge3$, a bound which is nontrivial when $y\ge x^{1-o(1)}$, and which starts to outclass \eqref{eq:P_{a,u,v}(x,y)} once $y\ge x\exp(-(\log x)^{1-\eta_0}(\log\log x)^{-1/2})$. In addition, the techniques used in \cite{Fan25b} may also be adapted to obtain nontrivial lower bounds for $P_{a,u,v}(x,y)$.

From Theorem \ref{thm:LMPshifted} one deduces at once the following result on linear forms $up+v$ in the image of the Carmichael function $\lambda(n)$, which may be defined as the largest order of a cyclic subgroup of the multiplicative group $(\Z/n\Z)^{\times}$. More explicitly, we have
 \[\lambda(p^{k})=\begin{cases}
\frac{1}{2}\varphi(p^k), & \text{if $p=2$ and $k\ge3$,}\\
\varphi(p^k), & \text{otherwise},
\end{cases}
\]
where $\varphi$ is the Euler totient function. The Chinese remainder theorem then yields $\lambda(n)=\text{lcm}\{\lambda(p^k)\colon p^k\parallel n\}$. We give the proof of this corollary right here since it is very easy.

\begin{cor}\label{cor:lambda}
Given any $u\in\N$ and $v\in\Z\setminus\{-u\}$, we have
\[\#((u\PP+v)\cap\lambda(\N)\cap[1,x])\ll_{u,v}\frac{\pi(x)}{((\log x)\log\log\log x)^{\eta_0}(\log\log x)^{1/2-\eta_0}}\]
for all $x\ge 16$. On the other hand, we have for sufficiently large $x$ that
\[\#((u\PP-u)\cap\lambda(\N)\cap[1,x])\asymp_{u}\pi(x).\]
\end{cor}
\begin{proof}
We start off by proving the first estimate. For each $n\in\N$, let $P^+(n)$ denote the largest prime factor of $n$, with the convention $P^+(1)=1$. The number of $n\le x$ with $P^+(n)\le y\colonequals x^{\log\log\log x/3\log\log x}$ is $O(\pi(x)/\log x)$, which is of a smaller order than the asserted bound. Suppose now that $P^+(up+v)>y$. If $up+v=\lambda(m)$, then $m$ must have a prime factor $q>y$, and thus $(q-1)\mid up+v$. By Theorem \ref{thm:LMPshifted}, the number of primes $p\le x$ such that $up+v$ has a shifted-prime divisor $q-1>y-1$ is
\[\ll\frac{\pi(x)}{(\log x)^{\eta_0}\sqrt{\log\log x}}\left(\frac{\log\log x}{\log\log\log x}\right)^{\eta_0}=\frac{\pi(x)}{((\log x)\log\log\log x)^{\eta_0}(\log\log x)^{1/2-\eta_0}},\]
which matches the first asserted estimate.

To prove the second estimate, which corresponds to the case $v=-u$, observe that the upper bound stems from counting $n\le x$ of the form $u(p-1)$. For the lower bound, note that any given $u\in\N$ divides some element of $\lambda(\N)$. For instance, $u$ divides $\lambda\left(\prod_{q^{\ell}\parallel u}q^{\ell+2}\right)$. Write $u=2^{k}\ell$, where $k\in\Z_{\ge0}$ and $\ell\in\N\setminus2\N$. If $\ell=1$, then $\lambda(2^{k+3}p)=u(p-1)$ for any $p\equiv 3\psmod{4}$, which produces the lower bound $(1/2+o(1))\pi(x)$. Suppose now that $\ell\ge3$. Since $\ell$ divides an element of $\lambda(\N)$, we may assume $\lambda(m)=\ell d$ for some fixed $d,m\in\N$. It is evident that $d\in2\N$. Fix $a\in(\Z/\ell\Z)^{\times}$ with $\gcd(1+ad,\ell)=1$. The number of primes $m<p\le x$ with $p\equiv 1+ad\psmod{\ell d}$ is asymptotically $(1+o(1))\pi(x)/\varphi(\ell d)$. In addition, since $p\nmid m$ and $\gcd(\lambda(m),\lambda(p))=\gcd(\ell d,p-1)=\gcd(\ell,a)d=d$, we have $\lambda(mp)=\lcm[\lambda(m),\lambda(p)]=\ell d(p-1)/d=\ell(p-1)$. For every $n\in\N$, denote by $v_2(n)$ the highest power of 2 dividing $n$. Then $1\le v_2(\lambda(m))\le v_2(\lambda(p))=v_2(p-1)$. Consequently,
\begin{align*}
\lambda\left(2^{k+v_2(\lambda(p))-v_2(m)+2}mp\right)&=\lcm\left[\lambda\left(2^{k+v_2(\lambda(p))-v_2(m)+2}m\right),\lambda(p)\right]\\
&=\lcm\left[\lambda\left(2^{k+v_2(\lambda(p))+2}\right),\lambda\left(m/2^{v_2(m)}\right),\lambda(p)\right]\\
&=\lcm\left[2^{k+v_2(\lambda(p))},\lcm\left[\lambda\left(m/2^{v_2(m)}\right),\lambda(p)\right]\right]\\
&=\lcm\left[2^{k+v_2(p-1)},\ell(p-1)\right]\\
&=u(p-1).
\end{align*}
Hence, there are at least $(1+o(1))\pi(x)/\varphi(\ell d)$ integers in $\lambda(\N)\cap[1,x]$ of the form $u(p-1)$, which verifies the lower bound for the case $\ell\ge3$.
\end{proof}

Corollary \ref{cor:lambda} provides a shifted-prime analogue to  
$\#(\lambda(\N)\cap[1,x])\le x/(\log x)^{\eta_0+o(1)}$ \cite[Theorem 1.1]{LP14}, which was later upgraded to an equality \cite[Theorem 1]{FLP14}. It is therefore reasonable to conjecture that given any $u\in\N$ and $v\in\Z\setminus\{-u\}$ sharing the same parity, we have $\#((u\PP+v)\cap\lambda(\N)\cap[1,x])=\pi(x)/(\log x)^{\eta_0+o(1)}$. 

Theorem \ref{thm:LMPshifted} also has interesting implications for the largest order $T_{CM}(d)$ of a torsion subgroup of an elliptic curve with complex multiplication over a degree $d$ number field. For example, one can show, by utilizing Theorem \ref{thm:LMPshifted} and adapting the proof of \cite[Theorem 1.7]{MPP17}, that 
\[\#\{p\le x\colon T_{CM}(up+v)>y\}\ll_{u,v}\frac{\pi(x)}{(\log y)^{\eta_0}\sqrt{\log\log y}}\]
for all $x,y\ge3$, where $u\in\N$ and $v\in\Z\setminus\{-u\}$. We invite the interested reader to verify this assertion as well as explore potential analogues to other results in \cite{MPP17} and beyond.

\medskip
\section{The Hardy--Ramanujan inequality for sifted sets: Proof of Theorem \ref{thm:HRS}}\label{S:HRS}
We dedicate this section to a proof of Theorem \ref{thm:HRS}. Our proof leverages the ideas from \cite{Tim95} concerning shifted primes. As explained in the introduction, it suffices to treat the case where $k\ge1$ and $M_f(x,E)\ge1/\beta$. In addition, we may assume that $x$ is sufficiently large. To begin with, we have by Stirling's formula that
\[\frac{t^k}{k!}\asymp\frac{t^ke^k}{k^{k+1/2}}\gg\frac{t^ke^{k/2}}{k^k}\]
for $t>0$. Elementary calculus shows that the right-hand side is strictly increasing for $k\le e^{-1/2}t$ and strictly decreasing for $k>e^{-1/2}t$. Taking $t=M_f(x,E)+O(1)\le A_1\log\log x+O(1)$, we have
\begin{equation}\label{eq:M_f(x,E)^k/k!}
\frac{(M_f(x,E)+O(1))^{k-1}}{(k-1)!}\gg\frac{1}{(\log x)^{O(1)}}
\end{equation}
for $k\le \beta M_{f}(x,E)$. This assures us that the bounds in Theorem \ref{thm:HRS} are at least $\gg x/(\log x)^{O(1)}$ and that any quantity of a smaller order can be safely ignored.

We start with the following lemma on the harmonic means of nonnegative multiplicative functions. Well-known results of this kind are often formulated for divisor-bounded multiplicative functions (e.g. \cite[Exercise 14.5\,(b)]{Kou19}). We give a proof for the sake of completeness.
\begin{lem}\label{lem:meanf(n)/n}
Let $x\ge2$ and $A,B>0$. If $f\colon\N\to\R_{\ge0}$ is a multiplicative function such that $f(p)\le A$ for all primes $p$ and
\[\sum_{p}\sum_{\ell\ge2}\frac{f(p^{\ell})}{p^{\ell}}\le B,\]
then
\[\sum_{n\le x}\frac{f(n)}{n}\asymp_{A,B}\exp\left(\sum_{p\le x}\frac{f(p)}{p}\right).\]
\end{lem}
\begin{proof}
The upper bound is easy:
\[\sum_{n\le x}\frac{f(n)}{n}\le\prod_{p\le x}\sum_{\ell\ge0}\frac{f(p^{\ell})}{p^{\ell}}\le\exp\left(\sum_{p\le x}\sum_{\ell\ge1}\frac{f(p^{\ell})}{p^{\ell}}\right)\ll\exp\left(\sum_{p\le x}\frac{f(p)}{p}\right).\]
To prove the lower bound, we define the multiplicative function $h$ via the convolution identity $A^{\omega}=f\ast h$. Then $h(p)=A-f(p)\in[0,A]$ for all primes $p$. Consequently, we have by the upper bound in the lemma that
\[\sum_{n\le x}\frac{\mu(n)^2A^{\omega(n)}}{n}\le\sum_{n\le x}\frac{f(n)}{n}\sum_{n\le x}\frac{\mu(n)^2h(n)}{n}\ll\exp\left(\sum_{p\le x}\frac{A-f(p)}{p}\right)\sum_{n\le x}\frac{f(n)}{n}.\]
Thus, the lower bound will follow readily if we have
\begin{equation}\label{eq:mu^2A^omega}
\sum_{n\le x}\frac{\mu(n)^2A^{\omega(n)}}{n}\gg(\log x)^{A}
\end{equation}
But this is a direct consequence of the Landau--Selberg--Delange method \cite[Theorem II.5.2]{Ten15} (or \cite[Theorem 13.2]{Kou19}) and partial summation. 

Due to the elementary nature of \eqref{eq:mu^2A^omega}, however, it may be interesting to seek an `elementary' proof. Firstly, it is not hard to verify
\begin{equation}\label{eq:mu^2f}
\sum_{n\le x}\frac{\mu(n)^2f(n)}{n}\asymp\sum_{n\le x}\frac{f(n)}{n}.
\end{equation} 
Indeed, we have
\[\sum_{m\le x}\frac{f(m)}{m}\le\sum_{\substack{dn\le x\\d\text{~squarefull}\\n\text{~squarefree}}}\frac{f(d)f(n)}{dn}\le \sum_{\substack{d\le x\\d\text{~squarefull}}}\frac{f(d)}{d}\sum_{n\le x}\frac{\mu(n)^2f(n)}{n}\ll\sum_{n\le x}\frac{\mu(n)^2f(n)}{n},\]
since
\[\sum_{\substack{d\le x\\d\text{~squarefull}}}\frac{f(d)}{d}\le\prod_{p}\left(1+\sum_{\ell\ge2}\frac{f(p^{\ell})}{p^{\ell}}\right)=\prod_{p}\left(1+O\left(\sum_{\ell\ge2}\frac{1}{p^{2\ell/3}}\right)\right)\ll1.\]
So we may add or remove the factor $\mu(n)^2$ at will. Next, let $\lfloor A\rfloor$ be the integer part of $A$, and denote by $\{A\}=A-\lfloor A\rfloor\in[0,1)$ its fractional part. Put $N=\lfloor A\rfloor+1\in\N$. The upper bound asserted in the lemma gives
\[\sum_{n\le x}\frac{(1-\{A\})^{\omega(n)}}{n}\ll (\log x)^{1-\{A\}},\]
from which it follows that
\begin{align*}
\sum_{n\le x}\frac{\mu(n)^2A^{\omega(n)}}{n}&\gg(\log x)^{\{A\}-1}\sum_{n\le x}\frac{\mu(n)^2A^{\omega(n)}}{n}\sum_{n\le x}\frac{(1-\{A\})^{\omega(n)}}{n}\\
&\ge(\log x)^{\{A\}-1}\sum_{n\le x}\frac{(\mu^2A^{\omega}\ast(1-\{A\})^{\omega})(n)}{n}\\
&\ge(\log x)^{\{A\}-1}\sum_{n\le x}\frac{\mu(n)^2N^{\omega(n)}}{n}\\
&\gg(\log x)^{\{A\}-1}\sum_{n\le x}\frac{N^{\omega(n)}}{n}.
\end{align*}
where the last inequality arises from an application of \eqref{eq:mu^2f}. Thus, it is sufficient to show
\begin{equation}\label{eq:N^omega}
\sum_{n\le x}\frac{N^{\omega(n)}}{n}\gg(\log x)^{N}.
\end{equation}
Observe that the left-hand side above is equal to
\[\sum_{n\le x}\frac{((N-1)^{\omega}\mu^2\ast1)(n)}{n}\ge\sum_{n\le \sqrt{x}}\frac{\mu(n)^2(N-1)^{\omega(n)}}{n}\sum_{n\le\sqrt{x}}\frac{1}{n}\gg\log x\sum_{n\le \sqrt{x}}\frac{(N-1)^{\omega(n)}}{n}.\]
Repeating this procedure $N-1$ times yields \eqref{eq:N^omega}. 
\end{proof}

We will also need the following Hardy--Ramanujan type inequality with harmonic weights.
\begin{lem}\label{lem:HRf(n)/n}
Let $x\ge2$ and $\alpha,C_1,C_2>0$. Let $E\subseteq\PP$ and $g\in\{\omega,\Omega\}$. Suppose that $f\colon\N\to\R_{\ge0}$ is a multiplicative function such that
\[\sum_{\substack{p\in E\\\ell\ge2}}\frac{f(p^{\ell})}{p^{\ell}}\le C_1\]
and
\[\sum_{\substack{p\notin E\\\ell\ge2}}\frac{f(p^{\ell})}{p^{\ell}}+1_{g=\Omega}\sum_{\substack{p\in E\\\ell\ge2}}\frac{f(p^{\ell})}{p^{\ell}}\alpha^{\ell}\le C_2.\]
Then 
\[\sum_{\substack{n\le x\\g(n,E)=k}}\frac{f(n)}{n}\ll_{C_2}\frac{(M_f(x,E)+C_1)^{k}}{k!}e^{M_f(x,E^c)}\]
for all $k\in\Z_{\ge0}$ when $g=\omega$ and for all integers $0\le k\le \alpha M_f(x,E)$ when $g=\Omega$. 
\end{lem}
\begin{proof}
For $g=\omega$, it is easy to see that
\[\sum_{\substack{n\le x\\\omega(n,E)=k}}\frac{f(n)}{n}\le\frac{1}{k!}\left(\sum_{\substack{p^{\ell}\le x\\p\in E}}\frac{f(p^{\ell})}{p^{\ell}}\right)^{k}\prod_{p\notin E}\sum_{\ell\ge0}\frac{f(p^{\ell})}{p^{\ell}}\ll \frac{(M_f(x,E)+C_1)^{k}}{k!}e^{M_f(x,E^c)},\]
as required. For $g=\Omega$, we have
\[\sum_{\substack{n\le x\\\Omega(n,E)=k}}\frac{f(n)}{n}\le\sum_{\substack{a\ge 1\text{~squarefull}\\\Omega(a,E)\le k}}\frac{f(a)}{a}\sum_{\substack{b\le x\\\omega(b,E)=k-\omega(a,E)}}\frac{\mu(b)^2f(b)}{b}.\]
Applying the lemma for the case $g=\omega$ with $\mu^2f$ in place of $f$, we see that 
\[\sum_{\substack{b\le x\\\omega(b,E)=k-\omega(a,E)}}\frac{\mu(b)^2f(b)}{b}\ll\frac{(M_f(x,E)+C_1)^{k-\Omega(a,E)}}{(k-\Omega(a,E))!}e^{M_f(x,E^c)}.\]
Since $k\le \alpha M_f(x,E)$ implies
\[\frac{1}{(k-\Omega(a,E))!}=\frac{k(k-1)\cdots(k-\Omega(a,E)+1)}{k!}\le\frac{(\alpha M_f(x,E))^{\Omega(a,E)}}{k!},\]
we have 
\begin{align*}
\sum_{\substack{n\le x\\\Omega(n,E)=k}}\frac{f(n)}{n}&\ll\frac{(M_f(x,E)+C_1)^{k}}{k!}e^{M_f(x,E^c)}\sum_{a\ge 1\text{~squarefull}}\frac{f(a)\alpha^{\Omega(a,E)}}{a}\\
&=\frac{(M_f(x,E)+C_1)^{k}}{k!}e^{M_f(x,E^c)}\prod_{\substack{p\in E}}\left(1+\sum_{\ell\ge2}\frac{f(p^{\ell})}{p^{\ell}}\alpha^{\ell}\right)\prod_{\substack{p\notin E}}\left(1+\sum_{\ell\ge2}\frac{f(p^{\ell})}{p^{\ell}}\right)\\
&\ll\frac{(M_f(x,E)+C_1)^{k}}{k!}e^{M_f(x,E^c)},
\end{align*}
proving the lemma for the case $g=\Omega$.
\end{proof}

Our next result is a corollary of \cite[Lemma 1]{Tim95} of the Bombieri--Vinogradov type, which may be of some independent interest.
\begin{lem}\label{lem:BVTim}
Let $x\ge2$, $A,A_1,C>0$, $B\ge0$, $A_2\colon \R_{>0}\to\R_{>0}$, and $\epsilon\in(0,1/2)$. Let $\{c_n\}_{n\le x}\subseteq\C$ be such that $|c_n|\le 1$ for $n\le x$, and let $f,h\in\mathscr{M}(A_1,A_2)$. For $U\ge e^{(\log x)^A}$ and $V,W>0$, put
\[\Delta\colonequals\sum_{d\le Q}h(d)\max_{y\le x}\max_{(a,d)=1}\Bigg|\sum_{\substack{pn\le y\\p>U,\,n>V,\,P^-(n)>W\\pn\equiv a\psmod{d}}}c_nf(n)(\log p)^{B}-\frac{1}{\varphi(d)}\sum_{\substack{pn\le y\\p>U,\,n>V\\P^-(n)> W}}c_nf(n)(\log p)^{B}\Bigg|\]
with $1\le Q\le x^{1/2-\epsilon}$, where $P^-(n)$ denotes the least prime factor of $n$, with the convention $P^-(1)=\infty$. Then 
\[\Delta\ll_{A,A_1,A_2,B,C,\epsilon}x(\log x)^{D_1}\left(\left(\frac{1}{\min\{U,V\}^{\frac{1}{4}}}+\frac{1}{(\log x)^{C}}\right)(\log x)^{D_2}+(\log x)^{D_3}\sqrt{\frac{\log Y}{X\log X}}\right),\]
where $X\colonequals \max\{2,\min\{U,W\}\}$, $Y=\max\{2,Q/X\}$, and 
\begin{align*}
D_1&\colonequals\frac{A_1(A_1+1)+B-1}{2},\\
D_2&\colonequals\frac{A_1^2+\max\{2B,1\}}{2}+1,\\
D_3&\colonequals\frac{A_1+\max\{B-1,0\}}{2}.
\end{align*}
\end{lem}
 \begin{proof}
This result is similar to \cite[Lemma 3]{Tim95} in essence and follows from \cite[Lemma 1]{Tim95} in a similar fashion. Suppose that $x$ is sufficiently large. We take $a_n=c_n1_{n>V,\,P^-(n)>W}f(n)$ and $b_p=1_{p>U}(\log p)^{B}$ in \cite[Lemma 1]{Tim95}, the latter of which fulfills \cite[Eq. (5), Lemma 2]{Tim95} by the classical Siegel--Walfisz type estimates for Dirichlet character sums \cite[Corollary 11.18]{MV07}. Invoking Theorem \ref{thmx:Poll20b} we obtain that
\[\sum_{n\le y}|a_n|^2\sum_{p\le x/y}b_p^2\ll x(\log x)^{A_1^2+\max\{2B-1,0\}-1}\]
for all $2\le y\le x$, and that
\[\frac{1}{x}\sum_{n\le 2x}|a_n|^4\sum_{p\le 2x/n}b_p^4\ll(\log x)^{\max\{4B-1,0\}}\sum_{n\le 2x}\frac{|a_n|^4}{n}\ll(\log x)^{A_1^4+\max\{4B-1,0\}}.\]
By \cite[Lemma 1]{Tim95} we find that
\[\Delta_1\colonequals\sum_{d\le Q}\max_{y\le x}\max_{(a,d)=1}\Bigg|\sum_{\substack{pn\le y\\pn\equiv a\psmod{d}}}a_nb_p-\frac{1}{\varphi(d)}\sum_{\substack{pn\le y\\(pn,d)=1}}a_nb_p\Bigg|\]
satisfies
\begin{equation}\label{eq:Delta1}
\Delta_1\ll\left(\sqrt{x/U}(\log x)^2+\sqrt{x/V}(\log x)^2+\sqrt{x}(\log x)^{2-2C}\right)\sqrt{x}(\log x)^{A_1^2/2+\max\{B,1/2\}}.
\end{equation}
Let 
\[\Delta_2\colonequals\sum_{d\le Q}\max_{y\le x}\max_{(a,d)=1}\Bigg|\sum_{\substack{pn\le y\\pn\equiv a\psmod{d}}}a_nb_p-\frac{1}{\varphi(d)}\sum_{pn\le y}a_nb_p\Bigg|.\]
Then $|\Delta_2-\Delta_1|$ is bounded above by
\[\sum_{d\le Q}\frac{1}{\varphi(d)}\sum_{q\mid d}\sum_{\substack{pn\le x\\q\mid pn}}|a_n|b_p\\
\le\sum_{d\le Q}\frac{1}{\varphi(d)}\left(\sum_{\substack{pn\le x\\p\mid d,\,p>U}}f(n)(\log p)^{B}+\sum_{\substack{q\mid d\\q>W}}\sum_{\substack{pn\le x\\q\mid n}}f(n)(\log p)^{B}\right).\]
By Theorem \ref{thmx:Poll20b}, the right-hand side is
\begin{align*}
&=\sum_{d\le Q}\frac{1}{\varphi(d)}\left(\sum_{\substack{p\mid d\\p>U}}(\log p)^{B}\sum_{n\le x/p}f(n)+\sum_{\substack{q\mid d\\q>W}}\sum_{\substack{n\le x\\q\mid n}}f(n)\sum_{p\le x/n}(\log p)^{B}\right)\\
&\ll x\sum_{d\le Q}\frac{1}{\varphi(d)}\left((\log x)^{A_1-1}\sum_{\substack{p\mid d\\p>U}}\frac{(\log p)^B}{p}+(\log x)^{\max\{B-1,0\}}\sum_{\substack{q\mid d\\q>W}}\sum_{\substack{n\le x\\q\mid n}}\frac{f(n)}{n}\right)\\
&\le x\sum_{d\le Q}\frac{1}{\varphi(d)}\left((\log x)^{A_1+B-1}\sum_{\substack{p\mid d\\p>U}}\frac{1}{p}+(\log x)^{\max\{B-1,0\}}\sum_{\substack{q\mid d\\q>W}}\sum_{q^{\ell}\le x}\frac{f(q^{\ell})}{q^{\ell}}\sum_{\substack{m\le x/q^{\ell}\\q\nmid m}}\frac{f(m)}{m}\right)\\
&\ll x(\log x)^{A_1+\max\{B-1,0\}}\sum_{d\le Q}\frac{1}{\varphi(d)}\left(\sum_{\substack{p\mid d\\p>U}}\frac{1}{p}+\sum_{\substack{q\mid d\\q>W}}\frac{1}{q}\right)\\
&\ll x(\log x)^{A_1+\max\{B-1,0\}}\sum_{\min\{U,W\}<p\le Q}\frac{1}{p}\sum_{\substack{d\le Q\\p\mid d}}\frac{1}{\varphi(d)}\\
&\le x(\log x)^{A_1+\max\{B-1,0\}}\sum_{\min\{U,W\}<p\le Q}\frac{1}{\varphi(p)p}\sum_{d\le Q/p}\frac{1}{\varphi(d)}\\
&\ll x(\log x)^{A_1+\max\{B-1,0\}}\frac{\log Y}{X\log X}.
\end{align*}
Hence,
\begin{equation}\label{eq:Delta2}
\Delta_2\ll|\Delta_1|+x(\log x)^{A_1+\max\{B-1,0\}}\frac{\log Y}{X\log X}.
\end{equation}
Finally, we pass from $\Delta_2$ to $\Delta$ via an application of the Cauchy--Schwarz inequality. Appealing to Brun--Tichmarsh and Shiu's theorem \cite[Theorem 1]{Shiu80}, we have
\begin{align*}
\sum_{\substack{pn\le x\\pn\equiv a\psmod{d}}}|a_n|b_p&\ll\sum_{\substack{n\le \sqrt{x}\\(n,d)=1}}|a_n|\sum_{\substack{p\le x/n\\pn\equiv a\psmod{d}}}(\log p)^B+\sum_{\substack{p\le \sqrt{x}\\p\nmid d}}(\log p)^B\sum_{\substack{n\le x/p\\pn\equiv a\psmod{d}}}f(n)\\
&\ll\frac{x(\log x)^{B-1}}{\varphi(d)}\sum_{n<x}\frac{|a_n|}{n}+\frac{x(\log x)^{A_1-1}}{\varphi(d)}\sum_{p\le \sqrt{x}}\frac{(\log p)^{B}}{p}\\
&\ll\frac{x(\log x)^{A_1+B-1}}{\varphi(d)}
\end{align*}
whenever $d\le Q$. In particular, the special case $d=1$ yields
\[\frac{1}{\varphi(d)}\sum_{pn\le x}|a_n|b_p\ll\frac{x(\log x)^{A_1+B-1}}{\varphi(d)}.\]
By Cauchy--Schwarz and Lemma \ref{lem:meanf(n)/n} applied to $h(n)^2n/\varphi(n)$, we see that
\begin{align*}
\Delta&\ll \sqrt{x}(\log x)^{(A_1+B-1)/2}\sum_{d\le Q}\frac{h(d)}{\sqrt{\varphi(d)}}\max_{y\le x}\max_{(a,d)=1}\Bigg|\sum_{\substack{pn\le y\\pn\equiv a\psmod{d}}}a_nb_p-\frac{1}{\varphi(d)}\sum_{pn\le y}a_nb_p\Bigg|^{1/2}\\
&\ll\sqrt{x}(\log x)^{(A_1+B-1)/2}\left(\sum_{d\le Q}\frac{h(d)^2}{\varphi(d)}\right)^{1/2}\sqrt{\Delta_2}\\
&\ll \sqrt{x\Delta_2}(\log x)^{(A_1^2+A_1+B-1)/2}.
\end{align*}
Combining this with \eqref{eq:Delta1} and \eqref{eq:Delta2} completes the proof of the lemma.
\end{proof}

Now we embark on the proof of Theorem \ref{thm:HRS}. Our target is
\begin{equation}\label{eq:f(n)sifted}
\sum_{\substack{n\in\Ss\\g(n,E)=k}}f(n).
\end{equation}
We shall frequently apply Theorem \ref{thmx:Poll20b} when sieving is needed. Let $P(t,E)$ denote the product of all primes $p\in E\cap[1,t]$, and set $P(t)\colonequals P(t,\PP)$. For $k=1$, we have
\[\sum_{\substack{n\in \Ss\\g(n)=1}}f(n)\le\sum_{\substack{p^{\ell}\le x\\p\le \sqrt{x}}}f(p^\ell)+\sum_{\substack{n\in S\\(n,P(\sqrt{x}))=1}}f(n).\]
Note that
\[\sum_{\substack{p^{\ell}\le x\\p\le \sqrt{x}}}f(p^\ell)\ll\sum_{\substack{p^{\ell}\le x\\p\le \sqrt{x}}}p^{\ell/4}\ll x^{1/4}\sum_{p\le\sqrt{x}}\frac{\log x}{\log p}\ll\frac{x^{3/4}}{\log x}.\]
In addition, Theorem \ref{thmx:Poll20b} applied to $f1_{(n,P(\sqrt{x}))=1}$ yields
\[\sum_{\substack{n\in \Ss\\(n,P(\sqrt{x}))=1}}f(n)\ll \frac{x}{\log x}\exp\left(\sum_{\sqrt{x}<p\le x}\frac{f(p)}{p}-\sum_{p\le x}\frac{\nu(p)}{p}\right)\ll\frac{x}{\log x}e^{-M_{\nu}(x)}.\]
Combining the two estimates above, we obtain
\[\sum_{\substack{n\in \Ss\\g(n)=1}}f(n)\ll\frac{x}{\log x}e^{-M_{\nu}(x)},\]
as desired. More generally, if $g(n,E)=1$, then $n=p^{\ell}m$, where $p\in E$, $\ell\in\N$, and $m\in\N$ with $p\nmid m$. Thus,
\[\sum_{\substack{n\in\Ss\\g(n,E)=1}}f(n)\le\sum_{\substack{p^{\ell}m\in\Ss\\p^{\ell}\le\sqrt{x},\,p\in E\\(m,E)=1}}f(p^\ell)f(m)+\sum_{\substack{p^{\ell}m\in\Ss\\p\le \sqrt{x},\,p\in E\\p^{\ell}>\sqrt{x}}}f(p^{\ell}m)+\sum_{\substack{n\in \Ss\\(n,P(\sqrt{x},E))=1}}f(n).\]
The middle sum is easily seen to be
\[\ll x^{1/8}\sum_{\substack{\sqrt{x}<d\le x\\d\text{~squarefull}}}\sum_{m<x/d}1\le x^{9/8}\sum_{\substack{d>\sqrt{x}\\d\text{~squarefull}}}\frac{1}{d}\ll x^{7/8}.\]
By Theorem \ref{thmx:Poll20b}, we have
\begin{align*}
\sum_{\substack{n\in \Ss\\(n,P(\sqrt{x},E))=1}}f(n)&\ll\frac{x}{\log x}\exp\left(\sum_{p\le x}\frac{f(p)1_{p\nmid P(\sqrt{x},E)}}{p}-\sum_{p\le x}\frac{\nu(p)}{p}\right)\\
&=\frac{x}{\log x}\exp\left(\sum_{\substack{p\le \sqrt{x}\\p\notin E}}\frac{f(p)}{p}+\sum_{\sqrt{x}<p\le x}\frac{f(p)}{p}-\sum_{p\le x}\frac{\nu(p)}{p}\right)\\
&\ll\frac{x}{\log x}e^{M_f(x,E^{c})-M_{\nu}(x)}.
\end{align*}
and
\begin{align*}
\sum_{\substack{p^{\ell}m\in\Ss\\p^{\ell}\le \sqrt{x},\,p\in E\\(m,E)=1}}f(p^\ell)f(m)&\ll\sum_{\substack{p^{\ell}\le \sqrt{x}\\p\in E}}\frac{f(p^\ell)}{p^{\ell}}\cdot\frac{x}{\log(x/p^{\ell})}\exp\left(\sum_{\substack{q\le x/p^\ell\\q\notin E}}\frac{f(q)}{q}-\sum_{\substack{q\le x/p^\ell\\q\ne p}}\frac{\nu(q)}{q}\right)\\
&\ll\frac{x}{\log x}e^{M_f(x,E^{c})-M_{\nu}(x)}\sum_{\substack{p^{\ell}\le x\\p\in E}}\frac{f(p^\ell)}{p^{\ell}}\\
&=\frac{x(M_f(x,E)+O(1))}{\log x}e^{M_f(x,E^{c})-M_{\nu}(x)}.
\end{align*}
Collecting the estimates above confirms the case $k=1$ of Theorem \ref{thm:HRS}. Here, the term $O(1)$ depends on $A_2$ only.

Suppose now that $k\ge2$. For convenience, we will replace $\Ss$ by $\Ss'=\Ss\setminus\Tt$, where 
\[\Tt\colonequals\left\{n\le x\colon n\text{~has a square factor~}d^2>y\text{~or~}n\text{~has a factor~}d>z\text{~with~}P^+(d)\le y\right\},\]
with $y\colonequals e^{\sqrt{\log x}}$ and $z\colonequals x^{1/8}$. This is legitimate, since the contribution from $n\in\Tt$ is of a smaller order compared to the upper bounds in Theorem \ref{thm:HRS}. To see this, we observe that
\begin{align*}
\sum_{a\ge1\text{~squarefull}}\frac{f(a)}{a^{3/4}}=\prod_{p}\left(1+\sum_{\ell\ge2}\frac{f(p^{\ell})}{p^{3\ell/4}}\right)&=\prod_{p}\left(1+O\left(\sum_{\ell\ge2}\frac{1}{p^{2\ell/3}}\right)\right)\\
&=\prod_{p}\left(1+O\left(\frac{1}{p^{4/3}}\right)\right)\ll1
\end{align*}
and that
\begin{align*}
\sum_{\substack{a\ge1\\P^+(a)\le y}}\frac{f(a)}{a^{1-\delta}}=\prod_{p\le y}\left(1+\frac{f(p)}{p^{1-\delta}}+\sum_{\ell\ge2}\frac{f(p^{\ell})}{p^{(1-\delta)\ell}}\right)&=\prod_{p\le y}\left(1+\frac{f(p)}{p^{1-\delta}}+O\left(\frac{1}{p^{4/3}}\right)\right)\\
&\ll \exp\left(\sum_{p\le y}\frac{f(p)}{p^{1-\delta}}\right)\\
&=\exp\left(\sum_{p\le y}\frac{f(p)}{p}\left(1+O\left(\delta\log p\right)\right)\right)\\
&\ll\exp\left(\sum_{p\le y}\frac{f(p)}{p}\right)\ll(\log y)^{A_1},
\end{align*}
where $\delta=1/\log y=1/\sqrt{\log x}$. Hence, we have by Theorem \ref{thmx:Poll20b} and Rankin's trick that
\begin{align*}
\sum_{\substack{n\le x\\d^2\mid n\\d>y}}f(n)\le\sum_{\substack{y<a\le x\\a\text{~squarefull}}}f(a)\sum_{\substack{b\le x/a\\b\text{~squarefree}}}f(b)&\ll x\sum_{\substack{y<a\le x\\a\text{~squarefull}}}\frac{f(a)}{a}\exp\left(\sum_{p\le x/a}\frac{f(p)-1}{p}\right)\\
&\ll x\sum_{\substack{y<a\le x\\a\text{~squarefull}}}\frac{f(a)}{a}\left(\log\frac{2x}{a}\right)^{A_1-1}\\
&\ll x(\log x)^{\max\{A_1-1,0\}}\sum_{\substack{y<a\le x\\a\text{~squarefull}}}\frac{f(a)}{a}\\
&\ll xy^{-1/4}(\log x)^{\max\{A_1-1,0\}}\sum_{a\ge1\text{~squarefull}}\frac{f(a)}{a^{3/4}}\\
&\ll xy^{-1/4}(\log x)^{\max\{A_1-1,0\}}
\end{align*}
and that
\begin{align*}
\sum_{\substack{n\le x\\d\mid n,\,d>z\\P^+(d)\le y}}f(n)=\sum_{\substack{z<a\le x\\P^+(a)\le y}}f(a)\sum_{\substack{b\le x/a\\P^-(b)>y}}f(b)&\ll x(\log x)^{\max\{A_1-1,0\}}\sum_{\substack{z<a\le x\\P^+(a)\le y}}\frac{f(a)}{a}\\
&\ll xz^{-\delta}(\log x)^{\max\{A_1-1,0\}}\sum_{\substack{a\ge1\\P^+(a)\le y}}\frac{f(a)}{a^{1-\delta}}\\
&\ll xz^{-\delta}(\log x)^{\max\{A_1-1,0\}}(\log y)^{A_1},
\end{align*}
both of which are negligible compared to the upper bounds in Theorem \ref{thm:HRS}. In what follows, we will focus on \eqref{eq:f(n)sifted} with $\Ss'$ in place of $\Ss$.

Let us first consider the contribution to \eqref{eq:f(n)sifted} from the terms with $n=dm\in\Ss'$, where $d\le\sqrt{x}$ is composed of primes $p\in E\cap[1,x^{1/3}]$, and $m$ is composed of primes $p\in (E\cap(x^{1/3},x])\cup E^c$. In particular, we have $\gcd(m,d)=1$ and $1_{m>1,\,E=\PP}\le g(m,E)<\log x/\log x^{1/3}=3$. The contribution to \eqref{eq:f(n)sifted} from these terms with $m=1$ is obviously at most
\[\sum_{d\le\sqrt{x}}f(d)\ll \sqrt{x}(\log x)^{A_1-1},\]
which is negligible. On the other hand, the contribution to \eqref{eq:f(n)sifted} from these terms with $m>1$ is at most
\[\sum_{j=1_{E=\PP}}^{2}\sum_{\substack{d\le \sqrt{x}\\p\mid d\Rightarrow p\in E\\g(d,E)=k-j}}f(d)\sum_{\substack{dm\in\Ss'\\(m,P(x^{1/3},E))=1}}f(m).\]
By Theorem \ref{thmx:Poll20b}, we have
\begin{align*}
\sum_{\substack{dm\in\Ss'\\(m,P(x^{1/3},E))=1}}f(m)&\ll\frac{x}{d\log(x/d)}\exp\left(\sum_{\substack{x^{1/3}<p\le x/d\\p\in E}}\frac{f(p)}{p}+\sum_{\substack{p\le x/d\\p\notin E}}\frac{f(p)}{p}-\sum_{\substack{p\le x/d\\p\nmid d}}\frac{\nu(p)}{p}\right)\\
&\ll \frac{x}{d\log x}e^{M_f(x,E^{c})-M_{\nu}(x)}\prod_{p\mid d}\left(1+\frac{\nu(p)}{p}\right).
\end{align*}
Summing this on $d$, we find that
\[\sum_{\substack{d\le\sqrt{x}\\p\mid d\Rightarrow p\in E\\g(d,E)=k-j}}f(d)\sum_{\substack{dm\in\Ss'\\(m,P(x^{1/3},E))=1}}f(m)\ll\frac{x}{\log x}e^{M_f(x,E^{c})-M_{\nu}(x)}\sum_{\substack{d\le x\\g(d,E)=k-j}}\frac{h(d)}{d},\]
where
\[h(d)\colonequals f(d)1_{p\mid d\Rightarrow p\in E}\prod_{p\mid d}\left(1+\frac{\nu(p)}{p}\right).\]
Invoking Lemma \ref{lem:HRf(n)/n}, we have
\[\sum_{\substack{d\le x\\g(d,E)=k-j}}\frac{h(d)}{d}\ll\frac{(M_f(x,E)+O(1))^{k-j}}{(k-j)!},\]
from which we deduce 
\[\sum_{\substack{d\le\sqrt{x}\\p\mid d\Rightarrow p\in E\\g(d,E)=k-j}}f(d)\sum_{\substack{dm\in\Ss'\\(m,P(x^{1/3},E))=1}}f(m)\ll\frac{x(M_f(x,E)+O(1))^{k-j}}{(k-j)!\log x}e^{M_f(x,E^{c})-M_{\nu}(x)}.\]
Summing this on $1_{E=\PP}\le j\le 2$, we see that the total contribution to \eqref{eq:f(n)sifted} from the terms in current consideration is
\[\ll \frac{x(M_f(x,E)+O(1))^{k-1_{E=\PP}}}{(k-1_{E=\PP})!\log x}e^{M_f(x,E^c)-M_{\nu}(x)},\]
matching the upper bounds in Theorem \ref{thm:HRS}. The term $O(1)$ depends at most on $A_1,A_2,v$.

It remains to estimate the contribution to \eqref{eq:f(n)sifted} from the terms with $n=dm\in\Ss'$, where $d>\sqrt{x}$ is composed of primes $p\in E\cap[1,x^{1/3}]$, and $m$ is composed of primes $p\in (E\cap(x^{1/3},x])\cup E^c$. We write $d=d_1d_2$ and $m=m_1m_2$, where $d_1,m_1$ satisfies $P^+(d_1),P^+(m_1)\le y$, and  $d_2,m_2$ are such that $P^-(d_2),P^-(m_2)>y$. Then $d_1,m_1\le z$, $d_2>\sqrt{x}/z$, and $d_2m_2$ is squarefree. Hence, the contribution to \eqref{eq:f(n)sifted} from the terms in current consideration is 
\begin{align*}
&\le\sideset{}{'}\sum_{\substack{d_1,m_1\le z\\}}f(d_1m_1)\sideset{}{'}\sum_{\substack{d_1m_1d_2m_2\in \Ss'\\g(d_1m_1d_2m_2,E)=k}}f(d_2m_2)\\
&\ll\frac{1}{\log x}\sideset{}{'}\sum_{\substack{d_1,m_1\le z\\}}f(d_1m_1)\sideset{}{'}\sum_{\substack{d_1m_1d_2m_2\in \Ss'\\g(d_1m_1d_2m_2,E)=k}}f(d_2m_2)\log d_2\\
&=\frac{1}{\log x}\sideset{}{'}\sum_{\substack{d_1,m_1\le z\\}}f(d_1m_1)\sum_{y<p\le x^{1/3}}\log p\sideset{}{'}\sum_{\substack{d_1m_1d_2m_2\in \Ss'\\p\mid d_2\\g(d_1m_1d_2m_2,E)=k}}f(d_2m_2),
\end{align*}
where the primed sums run over the specified variables $d_1,d_2,m_1,m_2$ satisfying the constraints described above. Peeling off a prime factor $y<p\le x^{1/3}$ of $d_2$ and putting $a=d_1m_1$ and $b=d_2m_2/p$, we have $a\le z^2$ with $P^+(a)\le y$, $b\ge d_2/p>x^{1/6}/z$, $P^-(b)>y$, $p\nmid b$, and $g(ab,E)=k-1$. Thus, the contribution to \eqref{eq:f(n)sifted} from the terms in current consideration is 
\begin{equation}\label{eq:sumf(a)f(b)log p}
\ll \frac{1}{\log x}\sideset{}{'}\sum_{\substack{a\le z^2\\}}f(a)\sum_{y<p\le x^{1/3}}\log p\sideset{}{'}\sum_{\substack{abp\in \Ss'\\g(ab,E)=k-1}}f(b).
\end{equation}
To estimate the inner double sum, we consider the multiset
\[\Rr_a\colonequals\left\{n=bp\le x/a\colon y<p\le x^{1/3},b>x^{1/6}/z, P^-(b)>y, g(ab,E)=k-1\right\},\]
where each $n$ appears as many times as it factors as $n=bp$, and define the multisequence $\{s_n\}$, supported on $\Rr_a$, by $s_n=f(b)\log p$ for $n\in\Rr_a$. Let 
\[X_a\colonequals\sum_{n\in\Rr_a}s_n,\]
and set
\[\Ww_a(d)\colonequals\left\{c\in\Z/d\Z\colon c\,\psmod{q}\in a^{-1}\Ee_q\text{~for every prime~}q\mid d\right\}\]
for squarefree $d\in\N$ with $\gcd(a,d)=1$ and $P^-(d)>2v+2$. In addition, we define $r_a(d)$ via
\[\sum_{\substack{n\in\Rr_a\\n\psmod{d}\in\Ww_a(d)}}s_n=\frac{\nu(d)}{\varphi(d)}X_a+r_a(d),\]
where we have extended naturally the definition of $\nu$ via
\[\nu(n)\colonequals\prod_{p\mid n}\nu(p)\]
for all $n\in\N$. Clearly, $\#\Ww_a(d)=\nu(d)$ and $\nu(d)/\varphi(d)<1/2$. By Selberg's sieve \cite[Theorem I.4.25]{Ten15} with trivial modifications to allow for the weights $s_n$, we have 
\[\sum_{y<p\le x^{1/3}}\log p\sideset{}{'}\sum_{\substack{abp\in \Ss'\\g(ab,E)=k-1}}f(b)\le\sum_{\substack{n\in\Rr_a\\q\nmid a,\, 2v+2<q\le w\Rightarrow n\notin\Ww_a(q)}}s_n\le X_aJ(D,w)^{-1}+\sum_{\substack{d\le D^2\\P^+(d)\le w\\(d,a)=1}}3^{\omega(d)}|r_a(d)|,\]
where we choose $w=D^2=(x/a)^{2/5}\ge (x/z^2)^{2/5}=x^{3/10}$, so that
\[J(D,w)\colonequals \sum_{\substack{n\le D^2\\P^+(n)\le w}}\frac{h_1(n)}{n}=\sum_{n\le (x/a)^{2/5}}\frac{h_1(n)}{n}\]
with
\[h_1(n)=\mu(n)^21_{(n,a)=1,\,P^-(n)>2v+2}\prod_{q\mid n}\frac{\nu(q)q}{q-1-\nu(q)}.\]
By Lemma \ref{lem:meanf(n)/n}, we have
\[J(D,w)\gg \exp\left(\sum_{\substack{2v+2<q\le w\\q\nmid a}}\frac{\nu(q)}{q-1-\nu(q)}\right)\gg h_2(a)^{-1}e^{M_{\nu}(x)},\]
where
\begin{equation}\label{eq:h_2(a)}
h_2(a)\colonequals\prod_{q\mid a}\left(1+\frac{\nu(q)}{q}\right).
\end{equation}
To bound the sum involving $r_a(d)$, note that
\[r_a(d)=\sum_{c\in\Ww_a(d)}\left(\sum_{\substack{n\in\Rr_a\\n\equiv c\psmod d}}s_n-\frac{1}{\varphi(d)}\sum_{n\in\Rr_a}s_n\right)\]
whenever $\gcd(d,a)=1$. It follows by Lemma \ref{lem:BVTim} with $B=1$, $U=W=y$, $V=x^{1/6}/z$, $Q=D^2=(x/a)^{2/5}$, $c_b=1_{g(ab,E)=k-1}$, $h(d)=3^{\omega(d)}\nu(d)$, and $C$ sufficiently large, that the sum involving $r_a(d)$ is bounded above by
\[\sum_{d\le Q}h(d)\max_{(c,d)=1}\Bigg|\sum_{\substack{pb\le x/a\\p>U,\,b>V,\,P^-(b)>W\\pb\equiv c\psmod{d}}}c_bf(b)\log p-\frac{1}{\varphi(d)}\sum_{\substack{pb\le x/a\\p>U,\,b>V\\P^-(b)> W}}c_bf(b)\log p\Bigg|\ll\frac{x}{a(\log x)^{G+A_1-1}}\]
for an arbitrary constant $G>0$. Therefore,
\[\sum_{y<p\le x^{1/3}}\log p\sideset{}{'}\sum_{\substack{abp\in \Ss'\\g(ab,E)=k-1}}f(b)\ll X_ah_2(a)\prod_{q\le x}\left(1-\frac{\nu(q)}{q}\right)+\frac{x}{a(\log x)^{G+A_1-1}}.\]
Inserting this into \eqref{eq:sumf(a)f(b)log p} and recalling the definition of $X_a$, we find that the contribution to \eqref{eq:f(n)sifted} from the terms in current consideration is 
\begin{align*}
&\ll \frac{1}{\log x}e^{-M_{\nu}(x)}\sum_{\substack{a\le z^2\\P^+(a)\le y}}f(a)h_2(a)\sum_{y<p\le x^{1/3}}\log p\sum_{\substack{abp\le x\\b>x^{1/6}/z, P^-(b)>y\\g(ab,E)=k-1}}f(b)+\frac{x}{(\log x)^{G+A_1}}\sum_{a\le z^2}\frac{f(a)}{a}\\
&\ll\frac{x}{\log x}e^{-M_{\nu}(x)}\sum_{\substack{a\le z^2\\P^+(a)\le y}}\frac{f(a)h_2(a)}{a}\sum_{\substack{b\le x/a\\P^-(b)> y\\g(ab,E)=k-1}}\frac{f(b)}{b}+\frac{x}{(\log x)^{G}}.
\end{align*}
Since $h_2(n)\ge1$ and every $n\in\N$ has a unique factorization $n=ab$ with $P^+(a)\le y<P^-(b)$, we infer by Lemma \ref{lem:HRf(n)/n} that
\[\sum_{\substack{a\le z^2\\P^+(a)\le y}}\frac{f(a)h_2(a)}{a}\sum_{\substack{b\le x/a\\P^-(b)> y\\g(ab,E)=k-1}}\frac{f(b)}{b}\le \sum_{\substack{n\le x\\g(n,E)=k-1}}\frac{f(n)h_2(n)}{n}\ll\frac{(M_f(x,E)+O(1))^{k-1}}{(k-1)!}e^{M_f(x,E^c)}.\]
Besides, in view of the discussion on \eqref{eq:M_f(x,E)^k/k!}, the term $x/(\log x)^G$ above is negligible as long as $G$ is large enough. Hence, We conclude that the contribution to \eqref{eq:f(n)sifted} from the terms in current consideration is 
\[\ll\frac{x(M_f(x,E)+O(1))^{k-1}}{(k-1)!\log x}e^{M_f(x,E^c)-M_{\nu}(x)},\]
which is acceptable. Again, the term $O(1)$ depends at most on $A_1,A_2,v$. This completes the proof of Theorem \ref{thm:HRS}.

\medskip
\section{Large deviations of $g(n,E)$: Proofs of Corollaries \ref{cor:WNO}--\ref{cor:SPNO}}\label{S:NOg(n,E)}
In this section we derive Corollaries \ref{cor:WNO}--\ref{cor:SPNO} from Theorem \ref{thm:HRS}. The tools developed along the way will also be useful in later sections. We start with the following generalization of \cite[Lemma 2.9]{MPP17} which also appears explicitly on \cite[p.\,7]{HT88}.
\begin{lem}\label{lem:siftOomega}
Assuming the hypotheses of Theorem \ref{thm:HRS}, we have
\[\sum_{n\in\Ss}f(n)z^{g(n,E)}\ll_{A_1,A_2,p_01_{g=\Omega},v,\beta}\frac{x}{\log x}e^{(z-1)M_f(x,E)+M_f(x)-M_{\nu}(x)}\]
for any $z\in(0,\beta]$.
\end{lem}
\begin{proof}
We may suppose that $x$ is sufficiently large. It is clear that the case $g=\omega$ follows immediately from Theorem \ref{thmx:Poll20b}. As much as we would also like to apply Theorem \ref{thmx:Poll20b} to $f(n)z^{\Omega(n,E)}$, we are not allowed to do so in that it does not necessarily satisfy Condition (ii) in the definition of $\mathscr{M}(A_1,A_2)$. To overcome this issue, we pass from $f(n)z^{\Omega(n,E)}$ to the more amenable function $\tau_{f,z}(n,E)$, which is defined via the formal Euler product identity
\[\sum_{n\ge1}\frac{\tau_{f,z}(n,E)}{n^s}=\prod_{p\in E}\left(1-p^{-s}\right)^{-zf(p)}\prod_{p\notin E}\sum_{\ell\ge0}\frac{f(p^{\ell})}{p^{\ell}}.\]
The special case $\tau_{1,z}(n,\PP)$ agrees with the classical $z$-fold divisor function $\tau_z(n)$. Applying the Taylor series expansion to each factor on the right-hand side of the identity above yields
\[\tau_{f,z}(p^{\ell},E)=1_{p\in E}\binom{w_p+\ell-1}{\ell}+1_{p\notin E}f(p^{\ell})\]
for all primes powers $p^{\ell}$, where $w_p\colonequals zf(p)\le \alpha_2 A_1$. It is easily seen by Stirling's formula that $\tau_{f,z}(n,E)$ does satisfy the conditions required by Theorem \ref{thmx:Poll20b}. Now, define the multiplicative function $h$ via the Dirichlet convolution identity $fz^{\Omega(\cdot,E)}=h\ast \tau_{f,z}(\cdot,E)$. Then $h$ is a multiplicative function with $h(p^{\ell})=0$ when $p\notin E$ and 
\[h(p^{\ell})=\sum_{j=0}^{\ell}\binom{-w_p+j-1}{j}f(p^{\ell-j})z^{\ell-j}=\sum_{j=0}^{\ell}(-1)^j\binom{w_p}{j}f(p^{\ell-j})z^{\ell-j}\]
when $p\in E$. It is clear that the terms corresponding to $j=0$ and $j=1$ cancel each other. In particular, $h(p)=0$ for all $p\in E$. Moreover, if $w_p\in[0,1)$, then 
\[\left|\binom{w_p}{j}\right|=\frac{w_p}{j}\cdot\frac{(1-w_p)\cdots(j-1-w_p)}{(j-1)!}\le\frac{w_p}{j}<\frac{1}{j}\]
for all $j\ge1$. Since $f(p^{\ell-j})\ll p^{\epsilon(\ell-j)}$ with
\[\epsilon\colonequals \min\left\{\frac{1}{2}\left(1-\frac{\log \alpha_2 }{\log p_0}\right),\frac{1}{3}\right\},\]
we have
\[|h(p^{\ell})|\ll \left(zp^{\epsilon}\right)^{\ell}+\sum_{j=1}^{\ell}\frac{\left(zp^{\epsilon}\right)^{\ell-j}}{j}\ll \max\left\{\left(zp^{\epsilon}\right)^{\ell},1+1_{zp^{\epsilon}=1}\log\ell\right\}\]
for all $p\in E$ and $\ell\ge2$. If $w_p\ge1$, then
\begin{align*}
\left|\binom{w_p}{j}\right|&=\frac{w_p(\lfloor w_p\rfloor+1-w_p)}{j(j-1)}\cdot\frac{(w_p-1)\cdots(w_p-\lfloor w_p\rfloor)}{\lfloor w_p\rfloor!}\cdot\frac{(\lfloor w_p\rfloor+2-w_p)\cdots(j-w_p-1)}{(\lfloor w_p\rfloor+1)\cdots(j-2)}\\
&\le\frac{w_p}{j(j-1)}
\end{align*}
for all $j\ge\lfloor w_p\rfloor+1$, since $w_p-1<\lfloor w_p\rfloor\le w_p$. Consequently, we have
\[|h(p^{\ell})|\ll \left(zp^{\epsilon}\right)^{\ell}+\sum_{j=\lfloor w_p\rfloor+1}^{\ell}\frac{\left(zp^{\epsilon}\right)^{\ell-j}}{j(j-1)}\ll \max\left\{\left(zp^{\epsilon}\right)^{\ell},1\right\}\]
for all $p\in E$ and $\ell\ge2$ in the case $w_p\ge1$. Combining the two cases, we conclude that
\begin{equation}\label{eq:h(p^l)}
|h(p^{\ell})|\ll\max\left\{\left(zp^{\epsilon}\right)^{\ell},1+1_{zp^{\epsilon}=1}\log\ell\right\}
\end{equation}
for all $p\in E$ and $\ell\ge2$.
	
With the introduction of $h$, we may write
\[\sum_{n\in\Ss}f(n)z^{\Omega(n,E)}\le\sum_{ab\in\Ss x}|h(a)|\tau_{f,z}(b,E),\]
the right-hand side of which is at most
\[\sum_{a\le\sqrt{x}}|h(a)|\sum_{ab\in\Ss}\tau_{f,z}(b,E)+\sum_{b\le\sqrt{x}}\tau_{f,z}(b,E)\sum_{a\le x/b}|h(a)|.\]
Since $\tau_{f,z}(p,E)=zf(p)1_{p\in E}+f(p)1_{p\notin E}$, we find by Theorem \ref{thmx:Poll20b} that 
\begin{align*}
\sum_{ab\in\Ss}\tau_{f,z}(b,E)&\ll\frac{x}{a\log(x/a)}\exp\left(\sum_{\substack{p\le x/a\\p\in E}}\frac{zf(p)}{p}+\sum_{\substack{p\le x/a\\p\notin E}}\frac{f(p)}{p}-\sum_{\substack{p\le x/a\\p\nmid a}}\frac{\nu(p)}{p}\right)\\
&\ll \frac{h_2(a)x}{a\log x}e^{(z-1)M_f(x,E)+M_f(x)-M_{\nu}(x)},
\end{align*}
where $h_2(a)$ is defined by \eqref{eq:h_2(a)}. Summing on $a\le\sqrt{x}$ and using \eqref{eq:h(p^l)}, we have
\begin{align*}
\sum_{a\le\sqrt{x}}|h(a)|\sum_{b\le x/a}\tau_{f,z}(b,E)&\ll\frac{x}{\log x}e^{(z-1)M_f(x,E)+M_f(x)-M_{\nu}(x)}\sum_{a\le\sqrt{x}}\frac{|h(a)|h_2(a)}{a}\\
&\le\frac{x}{\log x}e^{(z-1)M_f(x,E)+M_f(x)-M_{\nu}(x)}\prod_{p\in E}\left(1+O\left(\sum_{\ell\ge2}\frac{|h(p^{\ell})|}{p^{\ell}}\right)\right)\\
&\ll\frac{x}{\log x}e^{(z-1)M_f(x,E)+M_f(x)-M_{\nu}(x)},
\end{align*}
which coincides with the asserted bound. On the other hand, let us fix 
\[\max\left\{\frac{1}{2},\frac{\log\alpha_2 }{\log p_0}\right\}+\epsilon<\delta<1.\]
By Rankin's trick and \eqref{eq:h(p^l)}, we have
\begin{align*}
\sum_{a\le x/b}|h(a)|\le\left(\frac{x}{b}\right)^{\delta}\sum_{a\le x/b}\frac{|h(a)|}{a^{\delta}}\le\left(\frac{x}{b}\right)^{\delta}\prod_{p\in E}\left(1+\sum_{\ell\ge2}\frac{|h(p^{\ell})|}{p^{\delta\ell}}\right)\ll\left(\frac{x}{b}\right)^{\delta}.
\end{align*}
Summing on $b\le\sqrt{x}$ and applying Theorem \ref{thmx:Poll20b} and partial summation, we obtain
\[\sum_{b\le\sqrt{x}}\tau_{f,z}(b,E)\sum_{a\le x/b}|h(a)|\ll x^{\delta}\sum_{b\le \sqrt{x}}\frac{\tau_{f,z}(b,E)}{b^{\delta}}\ll x^{(1+\delta)/2}(\log x)^{\max\{\alpha_2 ,1\}A_1-1},\]
which is negligible in comparison to the asserted bound. Gathering the estimates above completes the proof of the case $g=\Omega$ of the lemma.
\end{proof}

\begin{lem}\label{lem:sifttail}
In addition to the hypotheses in Theorem \ref{thm:HRS}, suppose also that $\beta>1$ and that $M_f(x,E)\ge c_0>0$ whenever $x\ge x_0\ge2$. Then for all $x\ge x_0$ we have
\[\sum_{\substack{n\in\Ss\\g(n,E)\le(1-\delta) M_f(x,E)}}f(n)\ll_{A_1,A_2,c_0,v,x_0}\frac{x}{\log x}e^{M_f(x)-M_{\nu}(x)}\frac{e^{-Q(1-\delta)M_f(x,E)}}{\delta\sqrt{(1-\delta)M_f(x,E)}}\]
when $\delta\in(0,1)$ and
\[\sum_{\substack{n\in\Ss\\g(n,E)\ge(1+\delta) M_f(x,E)}}f(n)\ll_{A_1,A_2,c_0,p_01_{g=\Omega},v,x_0,\beta}\frac{x}{\log x}e^{M_f(x)-M_{\nu}(x)}\frac{e^{-Q(1+\delta)M_f(x,E)}}{\delta\sqrt{M_f(x,E)}}\]
when $\delta\in(0,\beta-1]$, where $Q$ is defined as in Corollary \ref{cor:ErdosMT}.
\end{lem}
\begin{proof}
The first estimate follows from Theorem \ref{thm:HRS} (with $p_0=2$ and $\alpha_1=\alpha_2=1$) and \cite[Lemma (4.5)]{Nor76}\footnote{Note the subtle difference between the definition of our $Q$ and that of $Q$ on \cite[p.\,686]{Nor76}.}:
\begin{align*}
\sum_{\substack{n\in\Ss\\g(n,E)\le (1-\delta) M_f(x,E)}}f(n)&\ll\frac{x}{\log x}e^{M_f(x)-M_{\nu}(x)}\sum_{0\le k\le(1-\delta)M_f(x,E)}\frac{M_f(x,E)^k}{k!}e^{-M_f(x,E)}\\
&\ll \frac{x}{\log x}e^{M_f(x)-M_{\nu}(x)}\frac{e^{-Q(1-\delta)M_f(x,E)}}{\delta\sqrt{(1-\delta)M_f(x,E)}}.
\end{align*}
To prove the second estimate, let $\beta_0=(\beta+c)/2$, where $c=2\alpha_11_{g=\omega}+p_01_{g=\Omega}$ . The sum to be estimated is at most
\[\sum_{\substack{n\in\Ss\\(1+\delta)M_f(x,E)\le g(n,E)\le\beta_0 M_f(x,E)}}f(n)+\sum_{\substack{n\in\Ss\\g(n,E)>\beta_0 M_f(x,E)}}f(n)\beta_0^{g(n,E)-\beta_0M_f(x,E)}.\]
By Theorem \ref{thm:HRS} and \cite[Lemma (4.7)]{Nor76}, we have 
\[\sum_{\substack{n\in\Ss\\(1+\delta)M_f(x,E)\le g(n,E)\le\beta_0 M_f(x,E)}}f(n)\ll\frac{x}{\log x}e^{M_f(x)-M_{\nu}(x)}\frac{e^{-Q(1+\delta)M_f(x,E)}}{\delta\sqrt{M_f(x,E)}}.\]
On the other hand, Lemma \ref{lem:siftOomega} implies that the second sum does not exceed
\[\beta_0^{-\beta_0M_f(x,E)}\sum_{n\in\Ss}f(n)\beta_0^{g(n,E)}\ll\frac{x}{\log x}e^{M_f(x)-M_{\nu}(x)}e^{-Q(\beta_0)M_f(x,E)}\]
The right-hand side is of a smaller order due to the relation $1+\delta\le\beta<\beta_0$, which results in
\[Q(\beta_0)-Q(1+\delta)\ge\int_{\beta}^{\beta_0}\log t\,dt\ge(\beta_0-\beta)\log\beta>0.\]
Putting the two bounds above together proves the second estimate.
\end{proof}

\begin{rmk}\label{rmk:sifttail}
Under certain circumstances one can show that
\begin{equation}\label{eq:sumf(n)}
	\sum_{n\in\Ss}f(n)\gg\frac{x}{\log x}e^{M_f(x)-M_{\nu}(x)}.
\end{equation}
If so, then Lemma \ref{lem:sifttail} implies that $g(n,E)$, when weighted by $f(n)$ over $\Ss$, has normal order $M_f(x,E)$, provided that $M_f(x,E)\to\infty$ as $x\to\infty$. More precisely, if $0<\lambda\le\sqrt{M_f(x,E)}/2$, then Lemma \ref{lem:sifttail} with $p_0=2$, $\alpha_1=\alpha_2=3/2$ and $\delta=\lambda/\sqrt{M_f(x,E)}$ yields
\[\left(\sum_{n\in\Ss}f(n)\right)^{-1}\sum_{\substack{n\in\Ss\\|g(n,E)-M_f(x,E)|\ge\lambda \sqrt{M_f(x,E)}}}f(n)\ll\lambda^{-1}\exp\left(-\frac{\lambda^2}{2}+O\left(\frac{\lambda^3}{\sqrt{M_f(x,E)}}\right)\right).\]
If $\lambda=\lambda(x)\to\infty$ and $\lambda\le\epsilon\sqrt{M_f(x,E)}$ for a sufficiently small $\epsilon>0$, then the right-hand side above approaches zero as $x\to\infty$. In particular, $g(n,E)$ has normal order $M_f(x,E)$.
\end{rmk}

\begin{proof}[Proof of Corollary \ref{cor:WNO}]
In view of Remark \ref{rmk:sifttail} above, it suffices to prove the lower bound \eqref{eq:sumf(n)}. Let $D$ and $\Pp$ be as in Corollary \ref{cor:WNO}. The fundamental lemma of sieve theory \cite[Theorem 19.1]{Kou19} provides a lower bound sieve weight function $\lambda^-\colon\N\to[-1,1]$, which is supported on squarefree positive integers $d\le D$ composed entirely of primes in $\Pp$, with the property that $\lambda^-(1)=1$, $(1\ast\lambda^-)(d)\le0$ for all $d>1$, and that
\begin{equation}\label{eq:lambda^-}
\sum_{d\mid\Pp}\frac{\lambda^-(d)h(d)}{d}\gg \prod_{p\in\Pp}\left(1-\frac{h(p)}{p}\right),
\end{equation}
where $h(d)\colonequals\nu(d)F(d)d/\varphi(d)\in[0,d)$ for each $d$ in the support of $\lambda^-$. Let
\[\Hh_d\colonequals\left\{a\in\Z/d\Z\colon a\,\psmod{p}\in \Ee_p\text{~for every prime~}p\mid d\right\}.\]
Then $\#\Hh_d=\nu(d)$. In addition, set $X\colonequals\sum_{n\le x}f(n)$, and write
\[\sum_{\substack{\substack{n\le x\\n\psmod{d}\in\Hh_d}}}f(n)=\frac{h(d)}{d}X+r_d\]
for each $d$ in the support of $\lambda^-$. Since
\[\sum_{d\mid \Pp}\lambda^-(d)1_{n\psmod{d}\in\Hh_d}\le 1_{n\psmod{p}\notin\Ee_p\text{~for all~}p\in\Pp}\]
for all $n\in\N$, it follows from \eqref{eq:fBV} and \eqref{eq:lambda^-} that
\begin{align*}
\sum_{n\in\Ss}f(n)=\sum_{n\le x}f(n)1_{n\psmod{p}\notin\Ee_p\text{~for all~}p\in\Pp}&\ge\sum_{n\le x}f(n)\sum_{d\mid \Pp}\lambda^-(d)1_{n\psmod{d}\in\Hh_d}\\
&=X\sum_{d\mid\Pp}\frac{\lambda^-(d)h(d)}{d}+\sum_{d\mid\Pp}\lambda^-(d)r_d\\
&\gg X\prod_{p\in\Pp}\left(1-\frac{h(p)}{p}\right)-\sum_{\substack{d\le D\\d\mid\Pp}}|r_d|\\
&\gg Xe^{-M_{\nu}(x)}-o\left(\frac{x}{\log x}e^{M_f(x)-M_{\nu}(x)}\right)
\end{align*}
by \eqref{eq:fBV}. But our assumption \eqref{eq:sumf(p)} and Lemma \ref{lem:meanf(n)/n} imply that
\begin{align*}
X\ge\sum_{m\le x^{1-\theta}}f(m)\sum_{x^{1-\theta}<p\le x/m}f(p)&=\sum_{m\le x^{1-\theta}}f(m)\left(\sum_{p\le x/m}f(p)+O\left(\pi\left(x^{1-\theta}\right)\right)\right)\\
&\gg \frac{x}{\log x}\sum_{m\le x^{1-\theta}}\frac{f(m)}{m}\\
	&\gg \frac{x}{\log x}e^{M_f(x)}.
\end{align*}
Hence, we have verified \eqref{eq:sumf(n)}, and therefore Corollary \ref{cor:WNO} as well.
\end{proof}

Next, we derive Corollary \ref{cor:QFNO} by putting these ideas into practice.

\begin{proof}[Proof of Corollary \ref{cor:QFNO}]
We start with the assertion about $\Aa_F$ in the corollary. We may suppose that $x$ is sufficiently large. Bernays \cite[pp.\,91--92]{Ber12} showed that
\begin{equation}\label{eq:A_F}
\#\Aa_F=(C_{\Delta}+o_{\Delta}(1))\frac{x}{\sqrt{\log x}}
\end{equation}
for any primitive, positive-definite binary quadratic form $F$ of discriminant $\Delta$, where $C_{\Delta}>0$ is a constant depending only on $\Delta$. 

We may suppose $\lambda>3\max\{1,1/\sqrt{M(p_0,E)}\}$, since the case $\lambda\le3\max\{1,1/\sqrt{M(p_0,E)}\}$ is clearly trivial. As before, each $n\in\N$ can be factored uniquely as $n=ms$, where $m$ and $s$ are the squarefull and squarefree part of $n$, respectively. For each $n\in\Aa_F$, there exist $X,Y\in\Z$ such that $n=F(X,Y)=aX^2+bXY+cY^2$. By the equations
\begin{align*}
4an&=(2aX+bY)^2-\Delta Y^2,\\
4cn&=(bX+2cY)^2-\Delta X^2,
\end{align*}
we see that if $p\nmid 2d\Delta$ is a prime factor of the squarefree part $s$ of $n$, then $(\Delta/p)=1$. Denote by $\Pp$ the set of primes $p$ satisfying $p\mid 2d\Delta$ or $(\Delta/p)=1$. Then $\gcd(s,\Pp^c)=1$. Since the Kronecker symbol $(\Delta/\cdot)$ is a non-principle character modulo $|\Delta|$, the series $\sum_{p}(\Delta/p)p^{-1}$ is convergent \cite[p.\,57]{Dav00}. Thus,
\[\sum_{\substack{p\le t\\p\in\Pp}}\frac{1}{p}=\frac{1}{2}\sum_{p\le t}\left(1+\left(\frac{\Delta}{p}\right)\right)\frac{1}{p}+O(1)=\frac{1}{2}\log\log t+O(1)\]
for $t\ge3$. It follows by Theorem \ref{thmx:Poll20b} with $f=\mu^21_{(s,\Pp^c)=1}$ that the number of $n\in\Aa_F$ whose squarefull parts exceed $e^{\lambda\sqrt{M(x,E)}}$ is at most
\begin{align*}
&\vspace*{-5mm}\sum_{\substack{e^{\lambda\sqrt{M(x,E)}}<m\le\sqrt{x}\\m\text{~squarefull}}}\sum_{s\le x/m}\mu(s)^21_{(s,\Pp^c)=1}+\sum_{\substack{\sqrt{x}<m\le x\\m\text{~squarefull}}}\sum_{s\le x/m}\mu(s)^2\\
&\ll \frac{x}{\sqrt{\log x}}\sum_{\substack{e^{\lambda\sqrt{M(x,E)}}<m\le\sqrt{x}\\m\text{~squarefull}}}\frac{1}{m}+x\sum_{\substack{\sqrt{x}<m\le x\\m\text{~squarefull}}}\frac{1}{m}\\
&\ll\frac{x}{\sqrt{\log x}}e^{-\lambda\sqrt{M(x,E)}/2},
\end{align*}
which is acceptable since $3<\lambda\le\sqrt{M(x,E)}/2$ and $M(x,E)\ge M(p_0,E)>0$. Hence, we may restrict ourselves to those $n\in\Aa_F$ whose squarefull parts are at most $e^{M(x,E)}$. Furthermore, the number of such $n=ms\in\Aa_F$ with $g(m,E)>\theta\lambda\sqrt{M(x,E)}$ is at most
\begin{align*}
\sum_{\substack{m\le e^{\lambda\sqrt{M(x,E)}}\\m\text{~squarefull}\\g(m,E)>\theta\lambda\sqrt{M(x,E)}}}\sum_{s\le x/m}\mu(s)^21_{(s,\Pp^c)=1}&\ll\frac{x}{\sqrt{\log x}}\sum_{\substack{m\ge1\text{~squarefull}\\g(m,E)>\theta\lambda\sqrt{M(x,E)}}}\frac{1}{m}\\
&\le \frac{x}{\sqrt{\log x}}\beta^{-\theta\lambda\sqrt{M(x,E)}}\sum_{\substack{m\ge1\\m\text{~squarefull}}}\frac{\beta^{g(m,E)}}{m}\\
&\ll\frac{x}{\sqrt{\log x}}\beta^{-\theta\lambda\sqrt{M(x,E)}},
\end{align*}
which is also acceptable as $\theta=1/(1+4\log\beta)$ fulfills the inequality $\theta\log\beta>(1-\theta)^2/4$. So we may focus on the rest of $n=ms\in\Aa_F$ with $m\le e^{\lambda\sqrt{M(x,E)}}$ and $g(m,E)\le \theta\lambda\sqrt{M(x,E)}$. For each such $n=ms$, the inequality $|g(n,E)-M(x,E)|\ge\lambda\sqrt{M(x,E)}$ implies $|g(s,E)-M(x,E)|\ge(1-\theta)\lambda\sqrt{M(x,E)}$. Since 
\[M(x,E)-M(x/m,E)\le\sum_{x/m<p\le x}\frac{1}{p}\ll\frac{\log m}{\log x}\le\frac{\lambda\sqrt{M(x,E)}}{\log x},\]
for all $m\le e^{\lambda\sqrt{M(x,E)}}$, we have
\begin{align*}
|g(s,E)-M(x/m,E)|&\ge\left(1-\theta+O\left(\frac{1}{\log x}\right)\right)\lambda\sqrt{M(x,E)}\\
&\ge\left(1-\theta+O\left(\frac{1}{\log x}\right)\right)\lambda\sqrt{M(x/m,E)}\,.
\end{align*}
Invoking Lemma \ref{lem:sifttail} with Remark \ref{rmk:sifttail} in mind, we find that the number of the remaining $n\in\Aa_F$ with $|g(n,E)-M(x,E)|\ge\lambda\sqrt{M(x,E)}$ is 
\begin{align*}
&\le\sum_{\substack{m\le e^{M(x,E)}\\m\text{~squarefull}}}\sum_{\substack{s\le x/m\\|g(s,E)-M(x/m,E)|\ge\eta\sqrt{M(x/m,E)}}}\mu(s)^21_{(s,\Pp^c)=1}\\
&\ll\frac{x}{\sqrt{\log x}}\sum_{\substack{m\le e^{M(x,E)}\\m\text{~squarefull}}}\frac{1}{m}\cdot\eta^{-1}\exp\left(-\frac{\eta^2}{2}+O\left(\frac{\eta^3}{\sqrt{M(x/m,E)}}\right)\right)\\
&\ll\frac{x}{\sqrt{\log x}}\cdot\eta^{-1}\exp\left(-\frac{\eta^2}{2}+O\left(\frac{\eta^3}{\sqrt{M(x,E)}}\right)\right)\\
&\ll\frac{x}{\sqrt{\log x}}\cdot\lambda^{-1}\exp\left(-\frac{((1-\theta)\lambda)^2}{2}+O\left(\frac{\lambda^3}{\sqrt{M(x,E)}}\right)\right),
\end{align*}
where $\eta=(1-\theta+O(1/\log x))\lambda$. Combining this with \eqref{eq:A_F} confirms the assertion about $\Aa_F$.

The proof of the assertion about $\Bb_F$ is similar. Instead of using \eqref{eq:A_F}, we appeal to \cite[Theorem 21.2]{FI10} to obtain $\#\Bb_F\asymp x/(\log x)^{3/2}$. In addition to the constraint $\gcd(s,\Pp^c)=1$, we also require that $n\not\equiv-k\psmod{p}$ for all $p\le\sqrt{x}$ with $p\nmid 2dkm\Delta$ when $n=ms>\sqrt{x}-k$. The rest of the proof runs in essentially the same way as above once we incorporate this additional information into our application of Theorem \ref{thmx:Poll20b} and Lemma \ref{lem:sifttail}, except that ``$x/\sqrt{\log x}$" gets replaced by $x/(\log x)^{3/2}$ and the weights $1/m$ in the summations become $1/\varphi(m)$.
\end{proof}

\begin{rmk}\label{rmk:QFNO}
For a reduced form $F(X,Y)=X^2+uY^2$, where $u=-\Delta/4\in\N$ is squarefree such that $\Q(\sqrt{-u})$ has class number 1, the characteristic function $f$ of numbers expressible by $F$ is multiplicative. Consequently, we see that Corollary \ref{cor:QFNO} holds with $\theta=0$ in this case. Using the same idea and \cite[Theorem 3]{Hoo74}, one can also show that Corollary \ref{cor:QFNO} holds for $\Cc_F\colonequals\{n\le x\colon \text{both~}n\text{~and~}n+k\text{~are sums of two squares}\}$ with any fixed $k\in\Z\setminus\{0\}$. 
\end{rmk}

We proceed to consider Corollary \ref{cor:SPNO}, whose proof rests on the following variant of \cite[Theorem 25.11]{FI10} on twin almost primes. For ease of notation, we write $P_2$ for a positive integer with at most two (not necessarily distinct) prime factors.
\begin{lem}\label{lem:TAP}
Let $A>0$, $x\ge2$ and $a\in\N\cap[1,(\log x)^A]$. There exist $x_0\ge2$, $\delta\in(0,1)$, and $k_0(x)\in\N\cap((\log x)^{3/2},\infty)$, all of which may depend on $A$, such that
\[\#\left\{p\le x\colon ap+b=cP_2\text{~with~}P^-(P_2)>x^{3/11}\right\}\ge\frac{W(abc)}{31c}\cdot\frac{C_2x}{(\log x)^2}\]
for all $x\ge x_0$ and all $b\in\Z\setminus\{0\}$ and $c\in\N$ such that $a|b|c\le x^{\delta}$, $2\mid abc$, $\gcd(a,b)=\gcd(c,ab)=1$, and $2c/\gcd(2,c)$ is indivisible by $k_0(x)$, where
\[W(n)\colonequals\prod_{\substack{p\mid n\\p>2}}\left(1-\frac{1}{p-1}\right)^{-1},\]
and
\[C_2\colonequals 2\prod_{p>2}\left(1-\frac{1}{(p-1)^2}\right)\]
is the twin prime constant.
\end{lem}
\begin{proof}
The lemma almost follows from the proof of \cite[Theorem 25.11]{FI10}. Given that only a few small changes need to be made, we elaborate on these changes rather than replicate the entire argument. Suppose that $x$ is sufficiently large and $\delta$ is sufficiently small. First of all, since we allow $a,b,c$ to vary with $x$, we use \cite[Theorem 22.5]{FI10} in place of \cite[Theorem 22.4]{FI10} when applying the linear sieve \cite[Theorems 12.19, 12.20]{FI10} with $r=1$, $D=x^{4/7-2c\delta}$ and $z=x^{3/11}$, where $c>0$ is an absolute constant. If $\delta\le3/11$, then our assumption $a|b|c\le x^{\delta}$ ensures that we still have
\[\prod_{\substack{2<p\le z\\p\nmid abc}}\left(1-\frac{1}{p-1}\right)\prod_{\substack{2<p\le z\\p\mid c}}\left(1-\frac{1}{p}\right)= V(z)W(abc)\gcd(2,c)\frac{\varphi(c)}{c}\]
in the treatment of $S(\Aa,z)$. The symmetric equation obtained by switching $a$ and $c$, which is needed in the estimation of $S(\Bb,z)$, also holds true. Next, the estimate for $|\Bb|$ is still valid thanks to the Siegel--Walfisz theorem on primes in arithmetic progressions \cite[Theorem 12.1]{Kou19}, since the modulus there is $2a/\gcd(2,c)\le 2(\log x)^A$ and $(ax+b)/c=(1+o(1))ax/c$. The last but most important change to be made concerns the estimation of 
\[|\Aa|=\pi\left(x;\frac{2c}{(2,a)},(c-b)\overline{a/(2,a)}\right),\]
where $(c-b)\overline{a/(2,a)}$ is the multiplicative inverse modulo $2c/(2,a)$. The modulus $2c/\gcd(2,a)$ is $\asymp x^{\delta}$ when $c\asymp x^{\delta}$, so Siegel--Walfisz is not applicable here. To surmount this hurdle, we invoke \cite[Proposition 8]{APR83} to find that for any given $\epsilon>0$, there exists $k_0(x)\in\N\cap((\log x)^{3/2},\infty)$ such that if $x$ is sufficiently large and $\delta$ is sufficiently small, then
\[|\Aa|>(1-\epsilon)\frac{x}{\varphi(2c/(2,a))\log x}\]
whenever $2c/\gcd(2,c)$ is not divisible by $k_0(x)$. This lower bound is sufficient for us to finish the proof as in that of \cite[Theorem 25.11]{FI10} once $\epsilon$ is taken to be small enough. 
\end{proof}

We are now in a position to carry out the proof of Corollary \ref{cor:SPNO}.
\begin{proof}[Proof of Corollary \ref{cor:SPNO}]
Suppose that $x$ is sufficiently large. Setting $n=ap+b$, we have
\[\sum_{\substack{p\le x\\|g(ap+b,E)-M_f(x,E)|\ge\lambda \sqrt{M_f(x,E)}}}f(ap+b)\ll\sum_{n\le a\sqrt{x}+b}f(n)+\sum_{\substack{a\sqrt{x}+b<n\le ax+b\\q\nmid a,\,q\le\sqrt{x}\Rightarrow n\not\equiv b\psmod{q}\\|g(n,E)-M_f(x,E)|\ge\lambda \sqrt{M_f(x,E)}}}f(n).\]
By Theorem \ref{thmx:Poll20b} and Lemma \ref{lem:sifttail} in conjunction with Remark \ref{rmk:sifttail}, the right-hand side is 
\[\ll\frac{x}{(\log x)^2}e^{M_f(x)}\cdot\lambda^{-1}\exp\left(-\frac{\lambda^2}{2}+O\left(\frac{\lambda^3}{\sqrt{M_f(x,E)}}\right)\right).\]
Hence, it suffices to show
\begin{equation}\label{eq:LBf(ap+b)}
\sum_{p\le x}f(ap+b)\gg\frac{x}{(\log x)^2}e^{M_f(x)}.
\end{equation}
This is where Lemma \ref{lem:TAP} comes into play; it tells us that there exist $\delta\in(0,3/11)$ and $k_0(x)\in\N\cap((\log x)^{3/2},\infty)$ such that
\[\#\left\{P_2\le (ax+b)/c\colon P^-(P_2)>x^{3/11}\text{~and~} ap+b=cP_2\text{~for some~}p\right\}\ge\frac{W(abc)}{31c}\cdot\frac{C_2x}{(\log x)^2}\]
for all $c\in\N\cap[1,x^{\delta}]$ with $2\mid abc$, $\gcd(c,ab)=1$, and $k_0(x)\nmid (2c/(2,c))$. Suppose first that $2\nmid ab$, so that $f(2^{\ell})\ge B$ for some $\ell\in\N$ by hypothesis. Then the left-hand side of \eqref{eq:LBf(ap+b)} is
\begin{align*}
&\ge\sum_{\substack{c\le x^{\delta}\\2^{\ell}\parallel c,\,k_0(x)\nmid c\\(c,ab)=1}}\mu\left(\frac{c}{2^{\ell}}\right)^2f(c)\sum_{\substack{P_2\le (ax+b)/c\\ap+b=cP_2\text{~for some~}p\\P^-(P_2)>x^{3/11}}}f(P_2)\\
&\ge\sum_{\substack{c\le x^{\delta}\\2^{\ell}\parallel c,\,k_0(x)\nmid c\\(c,ab)=1}}\mu\left(\frac{c}{2^{\ell}}\right)^2f(c)\left(\frac{W(abc)}{31c}\cdot\frac{B^2C_2x}{(\log x)^2}+O\left(\sum_{q\le\sqrt{(ax+b)/c}}\left(f(q^2)+1\right)\right)\right)\\
&\gg\frac{x}{(\log x)^2}\sum_{\substack{c\le x^{\delta}\\2^{\ell}\parallel c,\,k_0(x)\nmid c\\(c,ab)=1}}\mu\left(\frac{c}{2^{\ell}}\right)^2\frac{f(c)W(c)}{c}+O\left(\frac{\sqrt{x}}{\log x}\sum_{c\le x^{\delta}}\frac{f(c)}{\sqrt{c}}\right)\\
&\gg\frac{x}{(\log x)^2}\sum_{\substack{d\le x^{\delta}/2^{\ell}\\k_0(x)\nmid 2^{\ell}d\\(d,2ab)=1}}\mu(d)^2\frac{f(d)W(d)}{d}+O\left(x^{(1+\delta)/2}(\log x)^{A_1-2}\right).
\end{align*}
Let $r\colonequals k_0(x)/\gcd(k_0(x),2^{\ell})\asymp (\log x)^{3/2}$. Since Lemma \ref{lem:meanf(n)/n} implies that
\begin{align*}
\sum_{\substack{d\le x^{\delta}/2^{\ell}\\k_0(x)\mid 2^{\ell}d}}\mu(d)^2\frac{f(d)W(d)}{d}=\sum_{\substack{d\le x^{\delta}/2^{\ell}\\r\mid d}}\mu(d)^2\frac{f(d)W(d)}{d}&\le\mu(r)^2\frac{f(r)W(r)}{r}\sum_{d\le x^{\delta}/r}\mu(d)^2\frac{f(d)W(d)}{d}\\
&\ll\frac{f(r)W(r)}{r}e^{M_f(x)}\ll\frac{e^{M_f(x)}}{(\log x)^{4/3}},
\end{align*}
we have
\begin{align*}
\sum_{p\le x}f(ap+b)&\gg\frac{x}{(\log x)^2}\sum_{\substack{d\le x^{\delta}/2^{\ell}\\(d,2ab)=1}}\mu(d)^2\frac{f(d)W(d)}{d}+O\left(\frac{x}{(\log x)^{10/3}}e^{M_f(x)}\right)\\
&\gg\frac{x}{(\log x)^2}\exp\left(\sum_{\substack{2<p\le x^{\delta}/2^{\ell}\\p\nmid ab}}\frac{f(p)}{p}\right)\gg\frac{x}{(\log x)^2}e^{M_f(x)}
\end{align*}
by Lemma \ref{lem:meanf(n)/n} again. The case $2\mid ab$ is essentially the same: 
\begin{align*}
\sum_{p\le x}f(ap+b)&\ge\sum_{\substack{c\le x^{\delta}\\k_0(x)\nmid 2c\\(c,2ab)=1}}\mu(c)^2f(c)\sum_{\substack{P_2\le (ax+b)/c\\ap+b=cP_2\text{~for some~}p\\P^-(P_2)>x^{3/11}}}f(P_2)\\
&\gg\frac{x}{(\log x)^2}\sum_{\substack{c\le x^{\delta}\\k_0(x)\nmid 2c\\(c,2ab)=1}}\mu(c)^2\frac{f(c)W(c)}{c}+O\left(x^{(1+\delta)/2}(\log x)^{A_1-2}\right)\\
&\gg\frac{x}{(\log x)^2}\sum_{\substack{c\le x^{\delta}\\(c,2ab)=1}}\mu(c)^2\frac{f(c)W(c)}{c}+O\left(\frac{x}{(\log x)^{10/3}}e^{M_f(x)}\right)\\
&\gg\frac{x}{(\log x)^2}e^{M_f(x)}
\end{align*}
by the same calculations. This establishes \eqref{eq:LBf(ap+b)}, and therefore Corollary \ref{cor:SPNO}.
\end{proof}

Finally, we conclude this section by establishing Corollary \ref{cor:ErdosMT} as a consequence of Theorem \ref{thm:HRS} and Lemma \ref{lem:sifttail}.

\begin{proof}[Proof of Corollary \ref{cor:ErdosMT}]
Let $A=\max(A_1,1)$, and suppose that $x$ is sufficiently large. Observe that the contribution from those $n\in\Ss$ whose squarefull parts exceed $(\log x)^{3v+9A}$ is at most
\begin{align}\label{eq:sqf>(log x)^{3v+9A}}
\sum_{\substack{(\log x)^{3v+9A}<d\le x\\d\text{~squarefull}}}f(d)\sum_{\substack{m\le x/d\\m\text{~squarefree}}}f(m)\ll xe^{M_f(x)}\sum_{\substack{d>(\log x)^{3v+9A}\\d\text{~squarefull}}}\frac{f(d)}{d}&\ll \frac{x}{(\log x)^{v+3A}}e^{M_f(x)}\sum_{\substack{d\ge1\\d\text{~squarefull}}}\frac{f(d)}{d^{2/3}}\nonumber\\
&\ll\frac{x}{(\log x)^{v+3A}}e^{M_f(x)}.
\end{align}
Moreover, Lemma \ref{lem:sifttail} implies that the contribution from those $n\in\Ss$ with $\omega(n)>K\colonequals(1+\delta)M_f(x)$ with any $\delta\in(0,1)$ is 
\begin{equation}\label{eq:omega(n)>K}
\ll\frac{x}{\log x}e^{M_f(x)-M_{\nu}(x)}\frac{e^{-Q(1+\delta)M_f(x)}}{\delta\sqrt{M_f(x)}}.
\end{equation}
For each of the remaining $n=ab\in\Ss$ with $a,b\le\sqrt{x}$, let $d_1$ be the largest divisor of $a$ whose prime factors are also prime factors of $b$, and write $a=d_1m_1$. Symmetrically, we have the decomposition $b=d_2m_2$, where $d_2$ is the largest divisor of $b$ whose prime factors also divide $a$. Evidently, we have that $\gcd(d_1d_2,m_1m_2)=\gcd(m_1,m_2)=1$ and that $d_1d_2$ is squarefull. Since $\omega(n)\le K$ implies $\omega(m_1)+\omega(m_2)\le K$, it follows by Theorem \ref{thm:HRS} and \cite[Lemma (4.5)]{Nor76} that the contribution from the remaining $n=ab\in\Ss$ is 
\begin{align*}
&\le\sum_{\substack{d_1d_2\le(\log x)^{3v+9A}\\d_1d_2\text{~squarefull}}}f(d_1d_2)\sum_{\substack{k_1,k_2\ge0\\k_1+k_2\le K}}\sum_{\substack{m_1\le \sqrt{x}/d_1\\\omega(m_1)=k_1}}f(m_1)\sum_{\substack{d_1d_2m_1m_2\in\Ss\\m_2\le \sqrt{x}/d_2\\\omega(m_2)= k_2}}f(m_2)\\
&\ll\frac{\sqrt{x}}{\log x}e^{-M_{\nu}(x)}\sum_{\substack{d_1d_2\le(\log x)^{3v+9A}\\d_1d_2\text{~squarefull}}}\frac{f(d_1d_2)h_2(d_1d_2)}{d_2}\sum_{\substack{k_1,k_2\ge0\\k_1+k_2\le K}}\frac{M_f(x)^{k_2}}{k_2!}\sum_{\substack{m_1\le \sqrt{x}/d_1\\\omega(m_1)=k_1}}f(m_1)h_2(m_1)\\
&\ll\frac{x}{(\log x)^2}e^{-M_{\nu}(x)}\sum_{\substack{d_1,d_2\ge1\\d_1d_2\text{~squarefull}}}\frac{f(d_1d_2)h_2(d_1d_2)}{d_1d_2}\sum_{\substack{k_1,k_2\ge0\\k_1+k_2\le K}}\frac{M_f(x)^{k_1+k_2}}{k_1!k_2!}\\
&=\frac{x}{(\log x)^2}e^{2M_f(x)-M_{\nu}(x)}\sum_{\substack{d\ge1\\d\text{~squarefull}}}\frac{f(d)h_2(d)\tau(d)}{d}\sum_{0\le k\le K}\frac{(2M_f(x))^{k}}{k!}e^{-2M_f(x)}\\
&\ll\frac{x}{(\log x)^2}e^{2M_f(x)-M_{\nu}(x)}\frac{e^{-Q(c)\lambda}}{\sqrt{c(1-c)\lambda}}\\
&\ll\frac{x}{(\log x)^2}e^{2M_f(x)-M_{\nu}(x)}\frac{e^{-2Q(c)M_f(x)}}{\sqrt{(1-\delta)M_f(x)}},
\end{align*}
where $\lambda=2M_f(x)$, $c=(1+\delta)/2$, and $h_2$ is once again the multiplicative function defined by \eqref{eq:h_2(a)}. Combining this estimate with \eqref{eq:sqf>(log x)^{3v+9A}} and \eqref{eq:omega(n)>K}, we find that 
\begin{equation}\label{eq:sum_{n=ab}f(n)}
\sum_{n\in\Ss^{\ast}}f(n)\ll\left(\frac{e^{-Q(1+\delta)M_f(x)}}{\delta}+\frac{e^{(1-2Q(c))M_f(x)}}{(\log x)\sqrt{1-\delta}}\right)\frac{x}{(\log x)\sqrt{M_f(x)}}e^{M_f(x)-M_{\nu}(x)}.
\end{equation}
When $1/2<R_f(x)<1$, where we recall that $R_f(x)=(M_f(x)\log 2)/\log\log x$, we choose $\delta=1/R_f(x)-1\in(0,1)$ so that $(1-2Q(c))M_f(x)-\log\log x=-Q(1+\delta)M_f(x)$. Thus, \eqref{eq:sum_{n=ab}f(n)} becomes
\[\sum_{n\in\Ss^{\ast}}f(n)\ll\left(\frac{1}{1-R_f(x)}+\frac{1}{\sqrt{2R_f(x)-1}}\right)\frac{x}{(\log x)\sqrt{M_f(x)}}e^{(1-Q(1/R_f(x)))M_f(x)-M_{\nu}(x)}.\]
When $R_f(x)\le1/2$, we necessarily have $(1-2Q(c))M_f(x)-\log\log x<-Q(1+\delta)M_f(x)$ for any $\delta\in(0,1)$. Taking $\delta=1-c_0/(2M_f(x))$, we see that \eqref{eq:sum_{n=ab}f(n)} simplifies to
\[\sum_{n\in\Ss^{\ast}}f(n)\ll \left(1+\frac{4^{M_f(x)}\sqrt{M_f(x)}}{\log x}\right)\frac{x}{(\log x)\sqrt{M_f(x)}}e^{2(1-\log 2)M_f(x)-M_{\nu}(x)}.\]
This completes the proof of our corollary.
\end{proof}

\medskip
\section{The weighted normal order of $\omega(s(n))$: Proof of Theorem \ref{thm:EGPSTro}}\label{S:EGPSTro}
We devote this section to an analysis of $\omega(s(n))$. Despite being a bit overkill, the following Hardy--Ramanujan type inequality enables us to accomplish this goal. Plus, it seems possible to apply it to study $\omega(s_k(n))$ in a similar manner, where $s_k(n)\colonequals\sum_{d\mid n,d<n}d^k$ for $k\in\N$.
 \begin{prop}\label{prop:omega(Q(p))}
Let $Q_1,...,Q_r\in\Z[X]$ be distinct irreducible polynomials of positive degrees with $Q_j(0)\ne0$ for all $1\le j\le r$, such that $Q_0(X)\colonequals Q_1(X)\cdots Q_r(X)=\sum_{i=0}^{m}a_iX^i$ has degree $m$, discriminant $D$, height $\|Q\|\colonequals \max_{0\le i\le m}|a_i|$, and no fixed prime factors. Denote by $\rho_j(n)$ the number of solutions to $Q_j(X)\equiv0\psmod{n}$ for each $0\le j\le r$, and put
\[M_j\colonequals\sup_{t\ge2}\left|\sum_{p\le t}\frac{\rho_j(p)}{p}-\log\log t\right|\]
for each $1\le j\le r$. Let $\alpha,\delta\in(0,1)$ and $\beta>0$. Then there exist constants $c_0,K>0$, depending at most on $\alpha,\beta,\delta,m,r$, such that
\[\sum_{\substack{x-y<p\le x\\\forall j,\,\omega(Q_j(p))=k_j}}1\ll_{\alpha,\beta,\delta,m,r}\frac{|a_0|}{\varphi(|a_0|)}\left(\frac{|a_mD|}{\varphi_0(|a_mD|)}\right)^K\frac{e^My}{(\log x)^{r+1}}\prod_{j=1}^{r}\frac{(\log\log x+M_j)^{k_j-1}}{(k_j-1)!}\]
for all $x\ge c_0\|Q\|^{\delta}$, $x^{\alpha}<y\le x$, and $k_j\in\N\cap[1,\beta\log\log x]$, where $M\colonequals\sum_{j=1}^rM_j$ and
\[\varphi_0(n)\colonequals n\prod_{p\mid n}\left(1-\frac{\rho_0(p)}{p}\right).\]
\end{prop}
\begin{proof}
The proof follows exactly that of \cite[Th\'{e}or\`{e}me 1]{Ten18} with the index $p\in\PP$ in our sum replaced by $n\in\N$, at the cost of introducing the new constraint $q\mid n,q\nmid a_0\Rightarrow q>x^{\epsilon/3}$ to the sums $\sum^{\ast}$ and $\sum^{\ast\ast}$ on p.\,91 and p.\,93 of \cite{Ten18}, respectively. It is this additional constraint that produces the additional sieving factor
\[\prod_{\substack{q\le x^{\epsilon/3}\\q\nmid a_0}}\left(1-\frac{1}{q}\right)\ll\frac{|a_0|}{\varphi(|a_0|)\log x}\]
in the bounds for both sums, which in turn leads to the appearance of this factor in our bound compared to the bound in \cite[Th\'{e}or\`{e}me 1]{Ten18}. 
\end{proof}

The case $r=1$ of Proposition \ref{prop:omega(Q(p))} yields naturally the following corollary.
\begin{cor}\label{cor:omega(Q(p))}
Let $Q\in\Z[X]$ be irreducible of degree $m\in\N$ with leading coefficient $a_m$, constant coefficient $Q(0)=a_0\ne0$, discriminant $D$, height $\|Q\|$, and no fixed prime factors. Denote by $\rho_0(n)$ the number of solutions to $Q(X)\equiv0\psmod{n}$ and put
\[M\colonequals\sup_{t\ge2}\left|\sum_{p\le t}\frac{\rho_0(p)}{p}-\log\log t\right|.\]
For any $\alpha,\delta,\epsilon\in(0,1)$ and $\beta>1+\epsilon$, there exist constants $c_0,K>0$, depending at most on $\alpha,\beta,\delta,\epsilon,m$, such that
\[\sum_{\substack{x-y<p\le x\\|\omega(Q(p))-\log\log x|\ge\lambda\sqrt{\log\log x}}}1\ll_{\alpha,\beta,\delta,\epsilon,m}C_{\beta}\frac{|a_0|}{\varphi(|a_0|)}e^{\beta M}\cdot\frac{y}{\log x}\cdot\lambda^{-1}\exp\left(-\frac{\lambda^2}{2}+O_{\epsilon}\left(\frac{\lambda^3}{\sqrt{\log\log x}}\right)\right)\]
for all $x\ge c_0\|Q\|^{\delta}$, $x^{\alpha}<y\le x$ and $0<\lambda\le\epsilon\sqrt{\log\log x}$, where
\[C_{\beta}\colonequals\left(\frac{|a_mD|}{\varphi_0(|a_mD|)}\right)^Ke^{M}+e^{-M}\prod_{p\mid D}\left(1+\frac{1}{p}\right)^{m(m\beta+1)}\left(1+\beta\left(1+\frac{\rho_0(p)\beta}{p}\right)^{-1}\sum_{\ell\ge2}\frac{\rho_0(p^{\ell})}{p^{\ell}}\right)\]
with $\varphi_0$ defined as in Proposition \ref{prop:omega(Q(p))}.
 \end{cor}
 \begin{proof}
The proof is essentially the same as that of Lemma \ref{lem:sifttail} combined with Remark \ref{rmk:sifttail}, with the observation that
$(\log\log x+M)^{k-1}\le (\log\log x)^{k-1}e^{\beta M}$ for $1\le k\le\beta\log\log x$. Most importantly, one needs the estimate
\[\sum_{x-y<p\le x}\beta^{\omega(Q(p))}\ll\frac{|a_0|}{\varphi(|a_0|)}\cdot \frac{y}{\log x}\prod_{p\mid D}\left(1+\frac{1}{p}\right)^{m(m\beta+1)}\prod_{p\le x}\left(1-\frac{\rho_0(p)}{p}\right)\sum_{n\le x}\frac{\rho_0(n)\beta^{\omega(n)}}{n},\]
which is a consequence of \cite[Theorem 7 \& Remark of Corollary 1]{Hen12}. Since $\rho_0(p^{\ell})=\rho_0(p)\le m$ for $p\nmid D$ and $\ell\ge1$, we have
\begin{align*}
\sum_{n\le x}\frac{\rho_0(n)\beta^{\omega(n)}}{n}&\le\prod_{\substack{p\le x\\p\nmid D}}\left(1+\frac{\rho_0(p)\beta}{p-1}\right)\prod_{\substack{p\le x\\p\mid D}}\left(1+\beta\sum_{\ell\ge1}\frac{\rho_0(p^{\ell})}{p^{\ell}}\right)\\
&\ll\prod_{p\le x}\left(1+\frac{\rho_0(p)\beta}{p}\right)\prod_{p\mid D}\left(1+\beta\sum_{\ell\ge1}\frac{\rho_0(p^{\ell})}{p^{\ell}}\right)\left(1+\frac{\rho_0(p)\beta}{p}\right)^{-1}.
\end{align*}
Moreover, we have
\[\prod_{p\le x}\left(1-\frac{\rho_0(p)}{p}\right)\left(1+\frac{\rho_0(p)\beta}{p}\right)\le\prod_{p\le x}\left(1+(\beta-1)\frac{\rho_0(p)}{p}\right)\le e^{(\beta-1)M}(\log x)^{\beta-1}.\]
Piecing everything together gives the desired bound.
 \end{proof}
 
Undoubtedly, we also need to exploit to some extent the anatomy of $s(n)$ (and of $\sigma(n)$). In particular, we prove the following weighted analogues to \cite[Lemmas 2.2, 2.7]{Poll14}, respectively.

\begin{lem}\label{lem:p|sigma(n)}
Let $A_1>0$, $A_2\colon\R_{>0}\to\R_{>0}$, $\delta,\epsilon\in(0,1)$, and $f\in\mathscr{M}(A_1,A_2)$. Then 
\[\sum_{\substack{n\le x\\p\mid\sigma(n)}}f(n)\ll_{A_1,A_2,\delta,\epsilon}\left(\frac{1}{(\sqrt{p})^{1-\epsilon}}+\frac{\log\log x}{p}\right)\frac{x}{\log x}e^{M_f(x)}\]
for all $x\ge3$ and $p\in\PP\cap[2,x^{\delta}]$.
\end{lem}
\begin{proof}
We modify the proof of \cite[Lemma 2.2]{Poll14}. If $p\mid\sigma(n)$, then $p\mid\sigma(q^{\ell})$ for some prime power $q^{\ell}\parallel n$. Let $\theta=(1+\delta)/2$. If $\ell\ge2$, then $q^{\ell}>\sigma(q^{\ell})/2\ge p/2$. Thus, the contribution from the case $\ell\ge2$ is at most
\begin{align*}
\sum_{\substack{p/2<q^{\ell}\le x\\\ell\ge2}}\sum_{\substack{n\le x\\q^{\ell}\parallel n}}f(n)&=\sum_{\substack{p/2<q^{\ell}\le x\\\ell\ge2}}f(q^{\ell})\sum_{\substack{m\le x/q^{\ell}\\q\nmid m}}f(m)\\
&\le\sum_{\substack{p/2<q^{\ell}\le x^{\theta}\\\ell\ge2}}f(q^{\ell})\sum_{m\le x/q^{\ell}}f(m)+\sum_{\substack{x^{\theta}<q^{\ell}\le x\\\ell\ge2}}f(q^{\ell})\sum_{m\le x/q^{\ell}}f(m)\\
&\ll\frac{x}{\log x}e^{M_f(x)}\sum_{\substack{d>p/2\\d\text{~squarefull}}}\frac{1}{d^{1-\epsilon/2}}+xe^{M_f(x)}\sum_{\substack{d>x^{\theta}\\d\text{~squarefull}}}\frac{1}{d^{1-\epsilon/2}}\\
&\ll\frac{x}{(\sqrt{p})^{1-\epsilon}\log x}e^{M_f(x)}.
\end{align*}
On the other hand, if $\ell=1$, then $q\equiv-1\psmod{p}$. By Brun--Titchmarsh and partial summation, we find that the contribution from the case $\ell=1$ is at most
\begin{align*}
\sum_{\substack{q\le x\\q\equiv-1\psmod{p}}}\sum_{\substack{n\le x\\q\parallel n}}f(n)&=\sum_{\substack{q\le x\\q\equiv-1\psmod{p}}}f(q)\sum_{\substack{m\le x/q\\q\nmid m}}f(m)\\
&\ll xe^{M_f(x)}\sum_{\substack{q\le x\\q\equiv-1\psmod{p}}}\frac{1}{q\log(2x/q)}\\
&\ll\frac{x}{\varphi(p)}e^{M_f(x)}\left(\frac{1}{\log x}+\int_{p-1}^{x}\frac{dt}{t\log(2x/t)\log(3t/p)}\right)\\
&\ll\frac{x}{p}e^{M_f(x)}\left(\frac{1}{\log x}+\frac{1}{\log x}\int_{p-1}^{x^{\theta}}\frac{dt}{t\log(3t/p)}+\frac{1}{\log x}\int_{x^{\theta}}^{x}\frac{dt}{t\log(2x/t)}\right)\\
&\ll\frac{x\log\log x}{p\log x}e^{M_f(x)}.
\end{align*}
Combining the two estimates above finishes the proof of the lemma.
\end{proof}

\begin{lem}\label{lem:q|s(n)}
Let $A_1>0$, $A_2\colon\R_{>0}\to\R_{>0}$, and $f\in\mathscr{M}(A_1,A_2)$. Then
\[\sideset{}{'}\sum_{\substack{n\le x\\d\mid s(n)}}f(n)\ll_{A_1,A_2}\frac{x}{\varphi(d)\log(2y/z)}e^{M_f(x/y)}\sum_{d_1\mid d}d_1\left(\prod_{p\mid d_1}\sum_{\ell\ge v_{p}(d_1)}\frac{f(p^{\ell})}{p^{\ell}}\right)\]
for all $x\ge y\ge z\ge 1$ and $d\in\N\cap[1,z]$, where the index $n$ in the primed sum satisfies additionally $P^+(n)>y$ and $P^+(n)^2\nmid n$, and $v_p(d_1)$ is the highest power of $p$ dividing $d_1$.
\end{lem}
\begin{proof}
Write $n=mp$ with $p=P^+(n)$. Then $p>y$, $m\le x/y$, and $p\nmid m$, so that $s(n)=s(m)p+\sigma(m)$. Let $d_1=\gcd(d,s(m))$ and $d_2=d/d_1$. Then $d_1\mid\sigma(m)$, and so $d_1\mid m$. Thus,
\[\sideset{}{'}\sum_{\substack{n\le x\\d\mid s(n)}}f(n)\le\sum_{d_1d_1=d}\sum_{\substack{m\le x/y\\d_1\mid m}}f(m)\sum_{\substack{p\le x/m\\(s(m)/d_1)p\equiv-\sigma(m)/d_1\psmod{d_2}}}f(p),\]
which is
\begin{align*}
&\ll\sum_{d_1d_1=d}\sum_{\substack{m\le x/y\\d_1\mid m}}f(m)\cdot\frac{x}{m\varphi(d_2)\log(2x/(md_2))}\\
&\ll\frac{x}{\log(2y/z)}\sum_{d_1d_1=d}\frac{1}{\varphi(d_2)}\sum_{\substack{m_1\le x/y\\p\mid m_1\Rightarrow p\mid d_1}}\frac{f(m_1)}{m_1}\sum_{\substack{m_2\le x/m_1y\\(m_2,d_1)=1}}\frac{f(m_2)}{m_2}\\
&\ll\frac{x}{\log(2y/z)}e^{M_f(x/y)}\sum_{d_1d_1=d}\frac{1}{\varphi(d_2)}\prod_{p\mid d_1}\sum_{\ell\ge v_{p}(d_1)}\frac{f(p^{\ell})}{p^{\ell}}\\
&\le\frac{x}{\varphi(d)\log(2y/z)}e^{M_f(x/y)}\sum_{d_1\mid d}d_1\left(\prod_{p\mid d_1}\sum_{\ell\ge v_{p}(d_1)}\frac{f(p^{\ell})}{p^{\ell}}\right)
\end{align*}
by Brun--Titchmarsh.
\end{proof}

Finally, we derive from Lemma \ref{lem:p|sigma(n)} an upper bound for the mean value of $\omega(\gcd(\sigma(n),n))$.
\begin{cor}\label{cor:omega((sigma(n),n))}
Let $A_1>0$, $A_2\colon\R_{>0}\to\R_{>0}$, and $f\in\mathscr{M}(A_1,A_2)$. For sufficiently large $x$,
\[\sum_{n\le x}f(n)\omega((\sigma(n),n))\ll_{A_1,A_2}\frac{x}{\log x}e^{M_f(x)}\log_4x.\]
\end{cor}
\begin{proof}
Since $\omega((\sigma(n),n))\le\omega(n)\ll\log n$, the contribution from $n\le x$ whose squarefull parts exceed $x^{2/3}$ is at most
\begin{align*}
\log x\sum_{\substack{x^{2/3}<a\le x\\a\text{~squarefull}}}f(a)\sum_{\substack{b\le x/a\\b\text{~squarefree}}}f(b)&\ll xe^{M_f(x)}\log x\sum_{\substack{a>x^{2/3}\\a\text{~squarefull}}}\frac{f(a)}{a}\\
&\le x^{5/6}e^{M_f(x)}\log x\sum_{\substack{a>x^{2/3}\\a\text{~squarefull}}}\frac{f(a)}{a^{3/4}}\\
&\ll x^{5/6}e^{M_f(x)}\log x,
\end{align*}
which is negligible. So we only need to consider those $n\le x$ whose squarefull parts are at most $x^{2/3}$. The contribution from these $n$ does not exceed
\[\sum_{\substack{p^{\ell}\le x\\\ell\ge2\Rightarrow p^{\ell}\le x^{2/3}}}f(p^{\ell})\sum_{\substack{m\le x/p^{\ell}\\p\nmid m,\,p\mid\sigma(m)}}f(m)\le\sum_{p^{\ell}\le x^{2/3}}f(p^{\ell})\sum_{\substack{m\le x/p^{\ell}\\p\mid\sigma(m)}}f(m),\]
as $p\mid\sigma(m)$ implies $p\le\sigma(m)\le m^2\le (x/p)^2$, which then gives $p\le x^{2/3}$. Put $y=(\log\log x)^2$. By Lemma \ref{lem:p|sigma(n)} with $\delta=1/2$ and $\epsilon=1/3$, the right-hand side is
\[\ll\frac{x}{\log x}e^{M_f(x)}\left(\sum_{\substack{p^{\ell}\le x^{1/3}\\p\le y}}\frac{f(p^{\ell})}{p^{\ell}}+\sum_{\substack{p^{\ell}\le x^{1/3}\\p>y}}\frac{f(p^{\ell})}{p^{\ell+1/3}}+\sum_{x^{1/3}<p^{\ell}\le x^{2/3}}\frac{f(p^{\ell})}{p^{\ell}}\right)\ll \frac{x}{\log x}e^{M_f(x)}\log_4x,\]
as desired.
\end{proof}
 
We now have all the ingredients for proving Theorem \ref{thm:EGPSTro}. Assume $\lambda=c_0\sqrt{\log_4 x}$ with a constant $c_0>2$. Fix $\epsilon>0$ and suppose that $x$ is sufficiently large. Note that \eqref{eq:sumf(p)} implies that $M_f(x)\ge B\log\log x+O(1)$ and that
\[\sum_{n\le x}f(n)\gg\frac{x}{\log x}e^{M_f(x)}.\]
Let us fix a large $u>0$ and put $y=x^{1/u}$. Similarly to the proof of \cite[Theorem III.5.1]{Ten15}, Rankin's trick and the Halberstam--Richert theorem \cite[Theorem III.3.5]{Ten15} applied to $f(n)n^{\delta}1_{P^+(n)\le y}$ with $\delta=2/(3\log y)$ yield
\begin{align}\label{eq:P^+(n)>y}
\sum_{\substack{n\le x\\P^+(n)\le y}}f(n)&\le\sum_{n\le x^{3/4}}f(n)+ x^{-3\delta/4}\sum_{\substack{x^{3/4}<n\le x\\P^+(n)\le y}}f(n)n^{\delta}\nonumber\\
&\ll\frac{x^{3/4}}{\log x}e^{M_f(x)}+\frac{x^{1-3\delta/4}}{\log x}\sum_{\substack{n\le x\\P^+(n)\le y}}\frac{f(n)}{n^{1-\delta}}\nonumber\\
&\ll\frac{x^{3/4}}{\log x}e^{M_f(x)}+\frac{x^{1-3\delta/4}}{\log x}\exp\left(\sum_{p\le y}\frac{f(p)}{p^{1-\delta}}\right)\nonumber\\
&\ll\frac{x^{3/4}}{\log x}e^{M_f(x)}+\frac{x^{1-3\delta/4}}{\log x}\exp\left(\sum_{p\le y}\frac{f(p)}{p}+O\left(\delta\sum_{p\le y}\frac{\log p}{p}\right)\right)\nonumber\\
&\ll\left(x^{-1/4}+e^{-u/2}\right)\frac{x}{\log x}e^{M_f(x)},
\end{align}
where the implied constants are all independent of $u$. Moreover, we have
\begin{align}\label{eq:P^+(n)^2|n}
\sum_{\substack{n\le x\\P^+(n)>y\\P^+(n)^2\mid n}}f(n)\le\sum_{\substack{p^{\ell}\le x\\y<p\le\sqrt{x}\\\ell\ge2}}f(p^{\ell})\sum_{\substack{m\le x/p^{\ell}\\p\nmid m}}f(m)&\ll xe^{M_f(x)}\sum_{\substack{p^{\ell}\le x\\y<p\le\sqrt{x}\\\ell\ge2}}\frac{f(p^{\ell})}{p^{\ell}}\nonumber\\
&\ll xe^{M_f(x)}\sum_{p>y}\sum_{\ell\ge2}\frac{1}{p^{3\ell/4}}\ll \frac{xe^{M_f(x)}}{\sqrt{y}\log y}.
\end{align}
Besides, Corollary \ref{cor:omega((sigma(n),n))} implies that
\[\sum_{\substack{n\le x\\\omega((\sigma(n),n))>(\log_4 x)^2}}f(n)=O\left(\frac{x}{(\log x)\log_4 x}e^{M_f(x)}\right).\]
These contributions combined are bounded above by $(\epsilon/2)xe^{M_f(x)}/\log x$ when $u,x$ are large enough. Thus, it remains to consider the contribution from the set of $n\le x$ with $P^+(n)>y$, $P^+(n)^2\nmid n$, and $\omega(\gcd(\sigma(n),n))\le(\log_4 x)^2$. If we write $n=mp$ with $p=P^+(n)$, then $m\le x/y=x^{1-1/u}$ and $p\nmid m$, so that $s(n)=s(m)p+\sigma(m)$ and $d_m\colonequals\gcd(\sigma(m),s(m))=\gcd(\sigma(m),m)$ with $\omega(d_m)\le(\log_4 x)^2$. Let $a_m=s(m)/d_m$, $b_m=\sigma(m)/d_m$, and $Q_m(X)=a_mX+b_m$. Given each $m>1$, the contribution from this set of $n\le x$ is at most
\begin{equation}\label{eq:Q_m(p)}
f(m)\sum_{\substack{p\le x/m\\|\omega(Q_m(p))-\log\log x|\ge\lambda \sqrt{\log\log x}-\omega(d_m)}}f(p).
\end{equation}
Since $\log\log x-\log\log(x/m)=O(u)$ and $\omega(d_m)=o(\lambda \sqrt{\log\log x})$, the inequality $|\omega(Q_m(p))-\log\log x|\ge\lambda \sqrt{\log\log x}$ implies 
\[|\omega(Q_m(p))-\log\log(x/m)|\ge(1+o(1))\lambda\sqrt{\log\log x}\ge(1+o(1))\lambda\sqrt{\log\log(x/m)}\,.\]
By Corollary \ref{cor:omega(Q(p))} with $Q=Q_m$ and $\beta=c_0^2/4>1$ we see that \eqref{eq:Q_m(p)} is
\begin{align}\label{eq:f(m)a_mb_m}
&\ll \frac{f(m)x}{m\log(x/m)}\cdot\frac{b_m}{\varphi(b_m)}e^{(\beta+1) M}\cdot\lambda ^{-1}\exp\left(-(1+o(1))\frac{\lambda ^2}{2}+O\left(\frac{\lambda ^3}{\sqrt{\log\log(x/m)}}\right)\right)\nonumber\\
&\ll\frac{f(m)x}{m\log x}\cdot\frac{b_m}{\varphi(b_m)}\left(\frac{a_m}{\varphi(a_m)}\right)^{\beta+1}e^{-(1+o(1))\lambda^2/2},
\end{align}
provided that $x$ is sufficiently large, where
\[M=\sup_{t\ge2}\left|\sum_{p\le t}\frac{\rho_0(p)}{p}-\log\log t\right|\le\sum_{p\mid a_m}\frac{1}{p}+O(1).\]
This bound combined with \eqref{eq:P^+(n)>y}, \eqref{eq:P^+(n)^2|n}, and Lemmas \ref{lem:p|sigma(n)} and \ref{lem:q|s(n)}, shows, via the argument on \cite[p.\,144]{Poll14} with minor modifications, that for a suitably large constant $\alpha>0$ the sum of
\[\frac{f(m)}{m}\cdot\frac{b_m}{\varphi(b_m)}\left(\frac{a_m}{\varphi(a_m)}\right)^{\beta+1}\ll\frac{f(m)}{m}(\log\log x)^{\beta+2}\]
over those exceptional $1<m\le x$ failing 
\[\sum_{\substack{p\mid\sigma(m)\\p>(\log\log x)^{\alpha}}}\frac{1}{p}\le1\quad\text{~or}\quad\sum_{\substack{p\mid s(m)\\p>(\log\log x)^{\alpha}}}\frac{1}{p}\le1\] 
is $O(e^{M_f(x)})$. So the final contribution from these exceptional $m$ is 
\[\ll\frac{x}{\log x}e^{M_f(x)}\cdot e^{-(1+o(1))\lambda^2/2}.\]
Furthermore, since 
\[\frac{b_m}{\varphi(b_m)}\left(\frac{a_m}{\varphi(a_m)}\right)^{\beta+1}\ll\prod_{\substack{p\mid a_mb_m\\p\le(\log\log x)^{\alpha}}}\left(1+\frac{1}{p}\right)\prod_{\substack{p\mid a_m\\p\le(\log\log x)^{\alpha}}}\left(1+\frac{1}{p}\right)^{\beta}\ll(\log\log\log x)^{\beta+1}\]
for those non-exceptional $m>1$, the final contribution from them is
\[\ll\frac{x}{\log x}e^{M_f(x)}(\log\log\log x)^{\beta+1}e^{-(1+o(1))\lambda^2/2}.\]
Combining these contributions with the trivial contribution $O(x/\log x)$ from $m=1$, we find that the total remaining contribution is
\[\ll\left((\log\log\log x)^{\beta+1}e^{-(1+o(1))\lambda^2/2}+e^{-M_f(x)}\right)\frac{x}{\log x}e^{M_f(x)}<\frac{\epsilon}{2}\cdot\frac{x}{\log x}e^{M_f(x)}\]
for sufficiently large $x$, since $\beta+1<c_0^2/2$. This establishes Theorem \ref{thm:EGPSTro}.

\begin{rmk}\label{rmk:EGPSTro}
Perhaps the assumption $\lambda=c_0\sqrt{\log_4 x}$ can be weakened (ideally to $\lambda=\lambda(x)\to\infty$) by properly strengthening the hypotheses on $f$ and taking more advantage of the anatomy of $s(m)$ and $\sigma(m)$ when summing \eqref{eq:f(m)a_mb_m}.
\end{rmk}

\section{Shifted primes with a large shifted-prime divisor: Proof of Theorem \ref{thm:LMPshifted}}\label{S:SPlarge}
We end our paper with a proof of Theorem \ref{thm:LMPshifted}, which adapts that of \cite[Theorem 1.2]{MPP17} but is technically more involved. We begin with the following preparatory lemma on integers with exactly $k$ prime factors in an arithmetic progression.

\begin{lem}\label{lem:HRAP}
Let $A,\alpha_1>0$ and $\alpha_2\in(0,2)$. For $g\in\{\omega,\Omega\}$, put $\beta=\alpha_11_{g=\omega}+\alpha_21_{g=\Omega}$, and let $R\in(0,\beta]$. Then there exists a constant $c_0>0$, depending solely on $\beta$, such that
\begin{align}
\sum_{\substack{n\le x\\n\equiv a\psmod{d}\\g(n)=k}}1&\ll_{A,\beta}\frac{x(\log\log x)^{k-1}}{(k-1)!\varphi(d)\log x}+xe^{-c_0\sqrt{\log x}}\label{eq:HRAP1},\\
\sum_{\substack{n\le x\\n\equiv a\psmod{d}\\g(n)=k}}\frac{1}{n}&\ll_{A,R,\beta}\frac{(\log\log x)^k}{k!\varphi(d)}+(\log(3d))^RR^{-k},\label{eq:HRAP2}
\end{align}
uniformly for all $x\ge1$, $d\in\N\cap[1,e^{\sqrt{\log x}})$ with $\emph{rad}(d)\colonequals\prod_{p\mid d}p\le(\log x)^A$, $a\in\Z$ with $\gcd(a,d)=1$, and $k\in\N\cap[1,\beta\log\log x]$.
\end{lem}
\begin{proof}
We may assume $d$ is sufficiently large; otherwise the lemma follows immediately from the Hardy--Ramanujan inequality \eqref{eq:HR17} and Lemma \ref{lem:HRf(n)/n}. By the orthogonality relations for Dirichlet characters, we have
\[\sum_{\substack{n\le x\\n\equiv a\psmod{d}\\g(n)=k}}1=\frac{1}{\varphi(d)}\sum_{\chi\psmod{d}}\overline{\chi}(a)\sum_{\substack{n\le x\\g(n)=k}}\chi(n).\]
By \eqref{eq:HR17}, the contribution to the right-hand side from the principle character $\chi_0\psmod{d}$ is
\[\frac{1}{\varphi(d)}\sum_{\substack{n\le x\\g(n)=k}}\chi_0(n)\le\frac{1}{\varphi(d)}\sum_{\substack{n\le x\\g(n)=k}}1\ll\frac{x(\log\log x)^{k-1}}{(k-1)!\varphi(d)\log x}.\]
On the other hand, it can be shown \cite[Lemma 6]{Spiro85} that
\begin{equation}\label{eq:chi(n)z^{g(n)}}
\sum_{n\le x}\chi(n)z^{g(n)}\ll xe^{-c_1\sqrt{\log x}}
\end{equation}
for $\chi\ne\chi_0$ and $z\in\C$ with $|z|\le \beta$, where $c_1>0$ is a constant depending only on $\beta$. By Cauchy's integral formula we have
\[\sum_{\substack{n\le x\\g(n)=k}}\chi(n)=\frac{1}{2\pi i}\int_{|z|=\beta}\left(\sum_{n\le x}\chi(n)z^{g(n)}\right)z^{-k-1}\,dz\ll xe^{-c_1\sqrt{\log x}}\beta^{-k}\ll xe^{-c_0\sqrt{\log x}},\]
where $c_0\in(0,c_1)$ is a small constant depending on $\beta$. Collecting these estimates verifies \eqref{eq:HRAP1}.

The proof of \eqref{eq:HRAP2} is essentially the same. We make use of the identity
\begin{equation}\label{eq:HRAP3}
\sum_{\substack{n\le x\\n\equiv a\psmod{d}\\g(n)=k}}\frac{1}{n}=\frac{1}{\varphi(d)}\sum_{\chi\psmod{d}}\overline{\chi}(a)\sum_{\substack{n\le x\\g(n)=k}}\frac{\chi(n)}{n}.
\end{equation}
By Lemma \ref{lem:HRf(n)/n}, the contribution to the right-hand side from the principle character $\chi_0$ is
\begin{equation}\label{eq:HRAP4}
\frac{1}{\varphi(d)}\sum_{\substack{n\le x\\g(n)=k}}\frac{\chi_0(n)}{n}\le\frac{1}{\varphi(d)}\sum_{\substack{n\le x\\g(n)=k}}\frac{1}{n}\ll\frac{(\log\log x)^k}{k!\varphi(d)}.
\end{equation}
By \eqref{eq:chi(n)z^{g(n)}} and partial summation, we have
\[\sum_{n\le x}\frac{\chi(n)z^{g(n)}}{n}=\sum_{n\ge1}\frac{\chi(n)z^{g(n)}}{n}+O\left(e^{-c_0\sqrt{\log x}}\right)\]
for $\chi\ne\chi_0$ and $z\in\C$ with $|z|\le\beta$. It follows by Cauchy's integral formula that
\[\sum_{\substack{n\le x\\g(n)=k}}\frac{\chi(n)}{n}=\frac{1}{2\pi i}\int_{|z|=R}\left(\sum_{n\ge1}\frac{\chi(n)z^{g(n)}}{n}\right)z^{-k-1}\,dz+O\left(e^{-c_0\sqrt{\log x}}R^{-k}\right).\]
Since 
\[\sum_{n\ge1}\frac{\chi(n)z^{g(n)}}{n}=\prod_{p}\left(1+\sum_{\ell\ge1}\frac{\chi(p^{\ell})z^{g(p^{\ell})}}{p^{\ell}}\right)=\prod_{p}\left(1+\frac{\chi(p)z}{p}+O\left(\frac{1}{p^2}\right)\right)\ll |L(1,\chi)|^R,\]
and since $L(1,\chi)\ll\log d$ \cite[Lemma 11.2]{Kou19}, we have
\[\sum_{\substack{n\le x\\g(n)=k}}\frac{\chi(n)}{n}\ll(\log d)^RR^{-k}.\]
Inserting this and \eqref{eq:HRAP4} into \eqref{eq:HRAP3} completes the proof of \eqref{eq:HRAP2}.
\end{proof}

Next up is a technical result which sets the stage for the proof of Theorem \ref{thm:LMPshifted}.
\begin{lem}\label{lem:|an+b|/phi(|an+b|)}
Let $a,j\in\N$, $b\in\Z$, $\alpha_1>0$, and $\alpha_2\in(0,2)$. For $g\in\{\omega,\Omega\}$, put $\beta=\alpha_11_{g=\omega}+\alpha_21_{g=\Omega}$, and let $R\in(0,\beta]$. Then we have 
\begin{align}
\sum_{\substack{n\le x\\g(n)=k}}\left(\frac{|an+b|}{\varphi(|an+b|)}\right)^j&\ll_{a,b,j,\beta}\frac{x(\log\log x)^{k-1}}{(k-1)!\log x},\label{eq:|an+b|/phi(|an+b|)}\\
\sum_{\substack{n\le x\\g(n)=k}}\frac{1}{n}\left(\frac{|an+b|}{\varphi(|an+b|)}\right)^j&\ll_{a,b,j,\beta,R}\frac{(\log\log x)^{k}}{k!}+(\log\log x)^{R+j}R^{-k},\label{eq:1/phi(|an+b|)}
\end{align}
for all $x\ge3$ and $k\in\N\cap[1,\beta\log\log x]$.
\end{lem}
\begin{proof}
By the complete sub-multiplicativity of $n/\varphi(n)$, we may assume $\gcd(a,b)=1$. The case $b=0$ follows readily from Theorem \ref{thm:HRS} and Lemma \ref{lem:HRf(n)/n}. Suppose now that $b\ne0$. We first prove \eqref{eq:|an+b|/phi(|an+b|)}. The starting point is the inequality
\[\left(\frac{|an+b|}{\varphi(|an+b|)}\right)^j\ll\left(\frac{\sigma(|an+b|)}{|an+b|}\right)^j=\sum_{[\bdd]\mid(an+b)}\frac{1}{d_1\cdots d_j}\le\sum_{[\bdd]\mid(an+b)}\frac{1}{[\bdd]},\]
where $[\bdd]\colonequals\lcm[d_1,...,d_j]$ for the tuple $\bdd=(d_1,...,d_j)\in\N^j$. If we write $an+b=[\bdd]m$, then
\[\sum_{\substack{n\le x\\g(n)=k}}\left(\frac{|an+b|}{\varphi(|an+b|)}\right)^j\ll\sum_{\substack{\mathbf{d}\in\N^j\\ [\bdd]\le ax+|b|}}\frac{1}{[\bdd]}\sum_{\substack{m\le (ax+|b|)/[\bdd]\\ [\bdd]'m\equiv b'\psmod{a}\\g(([\bdd]m-b)/a)=k}}1\le \sum_{\substack{\mathbf{d}\in\N^j\\ [\bdd]\le ax+|b|}}\frac{1}{[\bdd]}\sum_{\ell=0}^{\min\{g(|b|),k\}}\sum_{\substack{m\le (ax+|b|)/[\bdd]\\ [\bdd]'m\equiv b'\psmod{a}\\g(([\bdd]'m-b')/a)=k-\ell}}1,\]
where $[\bdd]'=[\bdd]/\gcd([\bdd],b)$ and $b'=b/\gcd([\bdd],b)$. The contribution to the right-hand side from the terms with $\ell=k$ is obviously $O(1)$. Let $1\le c_{\bdd}<a$ be such that $[\bdd]'c_{\bdd}\equiv b'\psmod{a}$, and put $h_{\bdd}=([\bdd']c_{\bdd}-b')/a$. Then $h_{\bdd}\in\Z$ is coprime to $[\bdd]'$. Writing $m=am'+c_{\bdd}$ and $n'=[\bdd]'m'+h_{\bdd}$, we obtain
\begin{align*}
\sum_{\substack{n\le x\\g(n)=k}}\left(\frac{|an+b|}{\varphi(|an+b|)}\right)^j&\ll1+\sum_{\substack{\mathbf{d}\in\N^j\\ [\bdd]\le ax+|b|}}\frac{1}{[\bdd]}\sum_{\ell=0}^{\min\{g(|b|),k-1\}}\sum_{\substack{m'\ll x/[\bdd]\\g([\bdd]'m'+h_{\bdd})=k-\ell}}1\\
&\ll1+\sum_{\substack{\mathbf{d}\in\N^j\\ [\bdd]\le ax+|b|}}\frac{1}{[\bdd]}\sum_{\ell=0}^{\min\{g(|b|),k-1\}}\sum_{\substack{n'\ll x\\n'\equiv h_{\bdd}\psmod{[\bdd]'}\\g(n')=k-\ell}}1.
\end{align*}
Fix a large constant $A>0$. We have trivially
\begin{equation}\label{eq:lcm<=x}
\#\left\{\bdd\in\N^j\colon[\bdd]\le x\right\}=\sum_{n\le x}\#\left\{\bdd\in\N^j\colon[\bdd]=n\right\}\le\sum_{n\le x}\tau(n)^j\ll x(\log x)^{2^j-1}.
\end{equation}
By partial summation, we find that the contribution from the terms with $[\bdd]>(\log x)^{A}$ is
\[\ll x\sum_{\substack{\mathbf{d}\in\N^j\\ [\bdd]> (\log x)^{A}}}\frac{1}{[\bdd]^2}\ll x(\log x)^{2^j-1-A},\]
which is negligible compared to the right-hand side of \eqref{eq:|an+b|/phi(|an+b|)} when $A$ is sufficiently large in terms of $j$ and $\beta$. It now remains to estimate the contribution from the terms with $[\bdd]\le(\log x)^A$ and $\ell<k$. By \eqref{eq:HRAP1}, the contribution corresponding to each pair $(\bdd,\ell)$ is given by
\[\sum_{\substack{n'\ll x\\n'\equiv h_{\bdd}\psmod{[\bdd]'}\\g(n')=k-\ell}}1\ll\frac{x(\log\log x)^{k-1}}{(k-1)!\varphi([\bdd]')\log x}.\]
Summing this on $\bdd\in\N^j$ and $0\le\ell\le\min\{\omega(|b|),k-1\}$ and invoking \eqref{eq:lcm<=x}, we see that the remaining contribution is 
\[\ll\frac{x(\log\log x)^{k-1}}{(k-1)!\log x}\sum_{\bdd\in\N^j}\frac{1}{[\bdd]\varphi([\bdd]')}\ll\frac{x(\log\log x)^{k-1}}{(k-1)!\log x}\sum_{\bdd\in\N^j}\frac{1}{[\bdd]^{3/2}}\ll\frac{x(\log\log x)^{k-1}}{(k-1)!\log x},\]
which matches the desired upper bound. The proof of \eqref{eq:|an+b|/phi(|an+b|)} is now complete.

The proof of \eqref{eq:1/phi(|an+b|)} is similar. The same analysis shows that
\[\sum_{\substack{n\le x\\g(n)=k}}\frac{1}{n}\left(\frac{|an+b|}{\varphi(|an+b|)}\right)^j\ll\sum_{\substack{\mathbf{d}\in\N^j\\ [\bdd]\le ax+|b|}}\frac{1}{d_1\cdots d_j}\sum_{\ell=0}^{\min\{g(|b|),k\}}\sum_{\substack{n'\ll x\\n'\equiv h_{\bdd}\psmod{[\bdd]'}\\g(n')=k-\ell}}\frac{1}{n'}.\]
Since $d_1\cdots d_j\ge[\bdd]$, the contribution to the right-hand side from $[\bdd]>(\log x)^{A}$ is evidently $O((\log x)^{2^j-A})$ by \eqref{eq:lcm<=x}. Once $A$ is large enough in terms of $j$ and $\beta$, this contribution will be consumed by $(\log\log x)^k/k!$. To bound the remaining contribution, we apply \eqref{eq:HRAP2} (and its trivial case $k=0$) to get
\[\sum_{\substack{n'\ll x\\n'\equiv h_{\bdd}\psmod{[\bdd]'}\\g(n')=k-\ell}}\frac{1}{n'}\ll\frac{(\log\log x)^{k-\ell}}{(k-\ell)!\varphi([\bdd]')}+(\log(3[\bdd]'))^{R}R^{-k+\ell}.\]
Summing this on $[\bdd]\le (\log x)^{A}$ and $0\le\ell\le\min\{g(|b|),k\}$, we conclude that 
\begin{align*}
\sum_{\substack{n\le x\\g(n)=k}}\frac{1}{n}\left(\frac{|an+b|}{\varphi(|an+b|)}\right)^j&\ll\frac{(\log\log x)^{k}}{k!}\sum_{\mathbf{d}\in\N^j}\frac{1}{[\bdd]\varphi([\bdd]')}+(\log\log x)^{R}R^{-k}\sum_{\substack{\mathbf{d}\in\N^j\\ [\bdd]\le (\log x)^{A}}}\frac{1}{d_1\cdots d_j}\\
&\ll\frac{(\log\log x)^{k}}{k!}+(\log\log x)^{R}R^{-k}\sum_{\substack{\mathbf{d}\in\N^j\\ d_1\cdots d_j\le (\log x)^{jA}}}\frac{1}{d_1\cdots d_j}\\
&=\frac{(\log\log x)^{k}}{k!}+(\log\log x)^{R}R^{-k}\sum_{n\le (\log x)^{jA}}\frac{\tau_j(n)}{n}\\
&\ll\frac{(\log\log x)^{k}}{k!}+(\log\log x)^{R+j}R^{-k},
\end{align*}
which proves \eqref{eq:1/phi(|an+b|)}.
\end{proof}

We are now ready to prove Theorem \ref{thm:LMPshifted}. We may suppose throughout that $x,y$ are sufficiently large in terms of $a,u,v$. Besides, the letter $r$ will join $p$ and $q$ as notation for primes.

The case $v=0$ is easy. If $up=(q-a)m$, then $(q-a)/\gcd(q-a,u)$ divides $p$. So $q=\gcd(q-a,u)p+a$ if $q>u+|a|$. For each $d\mid u$, Brun's upper bound sieve shows that the number of $p\le x$ for which $dp+a\in\PP$ is $O(x/(\log x)^2)$. Summing this over all $d\mid u$, we arrive at a bound which is of a smaller order than our desired bound. Hence, we may assume $v\ne0$ from now on.

Suppose first that $x^{1/3}<y\le x$. Let $k=\lfloor (1/\log 2)\log\log x\rfloor$. By Lemma \ref{lem:sifttail} we have
\[\sum_{\substack{p\le x\\\Omega(up+v)>k}}1\ll\sqrt{x}+\sum_{\substack{u\sqrt{x}+v<n\le ux+v\\r\nmid uv,\,r\le\sqrt{x}\Rightarrow n\not\equiv v\psmod{r}\\\Omega(n)\ge(1/\log 2)\log\log x}}1\ll\frac{xe^{-Q(1/\log 2)(\log\log x+O(1))}}{(\log x)\sqrt{\log\log x}}\ll\frac{\pi(x)}{(\log y)^{\eta_0}\sqrt{\log\log y}},\]
which is the desired bound. Thus, we may focus on those $p\le x$ with $\Omega(up+v)\le k$. Write $up+v=(q-a)m$ with $q-a>y$, and suppose that $\omega(q-a)=i$ and $\omega(m)=j$ with $i+j\le k$. If $am+v=0$, then $up=mq$. For large $y$, this would imply that $p=q$ and $m=u$, so that $v=-am=-au$, contradicting our hypothesis. Hence, we always have $am+v\ne0$. Given $m\le (ux+v)/(q-a)\ll x^{2/3}$ and $0\le i\le k-j$, the count of such numbers $up+v$ is at most
\begin{align*}
\sum_{\substack{y<q-a\le (ux+v)/m\\(q-a)m\in u\PP+v\\\omega(q-a)=i}}1\le\sum_{\substack{y<n\le (ux+v)/m\\r\nmid a,\,r\le\sqrt{y}\Rightarrow n\not\equiv-a\psmod{r}\\r\nmid uvm,\,r\le\sqrt{y}\Rightarrow n\not\equiv \overline{m}v\psmod{r}\\\omega(n)=i}}1&\ll\frac{x}{m}\cdot\frac{(\log\log x)^{i}}{i!\log x}\prod_{\substack{r\le\sqrt{y}\\r\nmid auvm\\r\nmid(am+v)}}\left(1-\frac{2}{r}\right)\\
&\ll\left(\frac{m(am+v)}{\varphi(m)\varphi(|am+v|)}\right)^2\frac{x}{m}\cdot\frac{(\log\log x)^{i}}{i!(\log x)^3}
\end{align*}
by Theorem \ref{thm:HRS}. Summing this on $i,j\ge0$ with $i+j\le k$ and $m\le x$ with $\omega(m)=j\le k$, we find that the total count of such numbers $up+v$ is 
\[\ll\frac{x}{(\log x)^3}\sum_{\substack{i,j\ge0\\i+j\le k}}\frac{(\log\log x)^i}{i!}\sum_{\substack{m\le x\\\omega(m)=j}}\frac{1}{m}\left(\frac{m(am+v)}{\varphi(m)\varphi(|ma+v|)}\right)^2.\]
It follows from Cauchy--Schwarz, Lemma \ref{lem:HRf(n)/n}, and \eqref{eq:1/phi(|an+b|)} with $j=4$ and $R=1$, that 
\begin{align*}
\sum_{\substack{m\le x\\\omega(m)=j}}\frac{1}{m}\left(\frac{m(am+v)}{\varphi(m)\varphi(|ma+v|)}\right)^2&\ll\left(\sum_{\substack{m\le x\\\omega(m)=j}}\frac{1}{m}\left(\frac{m}{\varphi(m)}\right)^4\right)^{1/2}\left(\sum_{\substack{m\le x\\\omega(m)=j}}\frac{1}{m}\left(\frac{am+v}{\varphi(|am+v|)}\right)^4\right)^{1/2}\\
&\ll\frac{(\log\log x)^{j}}{j!}+(\log\log x)^{5/2}\left(\frac{(\log\log x)^{j}}{j!}\right)^{1/2}.
\end{align*}
Therefore, the total count of such numbers $up+v$ is 
\begin{align*}
&\ll\frac{x}{(\log x)^3}\sum_{\substack{i,j\ge0\\i+j\le k}}\frac{(\log\log x)^{i+j}}{i!j!}+\frac{x(\log\log x)^{5/2}}{(\log x)^3}\sum_{i=0}^{k}\frac{(\log\log x)^{i}}{i!}\sum_{j=0}^{k}\left(\frac{(\log\log x)^{j}}{j!}\right)^{1/2}\\
&\ll\frac{x}{(\log x)^3}\sum_{\ell=0}^{k}\frac{(2\log\log x)^{\ell}}{\ell!}+\frac{x\sqrt{k}(\log\log x)^{5/2}}{(\log x)^3}\sum_{i=0}^{k}\frac{(\log\log x)^{i}}{i!}\left(\sum_{j=0}^{k}\frac{(\log\log x)^{j}}{j!}\right)^{1/2}\\
&\ll\frac{x}{\log x}\cdot\frac{e^{-2Q(1/(2\log 2))\log\log x}}{\sqrt{\log\log x}}+\frac{x(\log\log x)^{3}}{(\log x)^3}\cdot(\log x)^{3/2}\\
&\ll\frac{\pi(x)}{(\log y)^{\eta_0}\sqrt{\log\log y}}
\end{align*}
by Cauchy--Schwarz and \cite[Lemma (4.5)]{Nor76}. This establishes Theorem \ref{thm:LMPshifted} for $x^{1/3}<y\le x$.

Next, we show that for $3\le y\le x^{1/3}$ one has
\begin{equation}\label{eq:y<q-a<=y^2}
P_{a,u,v}(x,y)-P_{a,u,v}(x,y^2)\ll\frac{\pi(x)}{(\log x)^{\eta_0}\sqrt{\log\log y}}.
\end{equation}
Let $\ell=\lfloor (1/\log 2)\log\log y\rfloor$, and define $\Omega_y(n)\colonequals\Omega(n,\PP\cap[2,y])$. By Lemma \ref{lem:sifttail} we have
\[\sum_{\substack{p\le x\\\Omega_y(up+v)>\ell}}1\ll\sqrt{x}+\sum_{\substack{u\sqrt{x}+v<n\le ux+v\\r\nmid uv,\,r\le\sqrt{x}\Rightarrow n\not\equiv v\psmod{r}\\\Omega_y(n)\ge(1/\log 2)\log\log y}}1\ll\frac{xe^{-Q(1/\log 2)(\log\log y+O(1))}}{(\log x)\sqrt{\log\log y}}\ll\frac{\pi(x)}{(\log y)^{\eta_0}\sqrt{\log\log y}},\]
which is again the desired bound. So we need only consider the contribution to \eqref{eq:y<q-a<=y^2} from those $p\le x$ for which $\Omega_y(up+v)\le \ell$. Write $up+v=(q-a)m$ with $y<q-a\le y^2$, and suppose that $\Omega_y(q-a)=i$ and $\Omega_y(m)=j$ with $i+j\le \ell$. Then $x/(q-a)\ge x/y^2\ge x^{1/3}\ge y$. Given $y<q-a<y^2$ and $0\le j\le \ell-i$, the corresponding contribution is at most
\begin{align*}
\sum_{\substack{m\le (ux+v)/(q-a)\\(q-a)m\in u\PP+v\\\Omega_y(m)=j}}1&\ll x^{1/4}+\sum_{\substack{u\sqrt{x}+v<m\le (ux+v)/(q-a)\\r\nmid uv(q-a),\,r\le x^{1/4}\Rightarrow m\not\equiv\overline{(q-a)}v\psmod{r}\\\Omega_y(m)=j}}1\\
&\ll x^{1/4}+\frac{x(\log\log y+O(1))^{j}}{j!q\log x}\exp\left(\sum_{y<r\ll x/(q-a)}\frac{1}{r}-\sum_{\substack{r\le x^{1/4}\\r\nmid uv(q-a)}}\frac{1}{r}\right)\\
&\ll\frac{|q-a|}{\varphi(|q-a|)q}\cdot\frac{\pi(x)(\log\log y)^{j}}{j!\log y}
\end{align*}
by Theorem \ref{thm:HRS}. Summing this on $i,j\ge0$ with $i+j\le \ell$ and $y<q-a\le y^2$ with $\Omega_y(q-a)=i$, we find that the contribution to \eqref{eq:y<q-a<=y^2} from those $p\le x$ for which $\Omega_y(up+v)\le \ell$ is 
\[\ll\frac{\pi(x)}{\log y}\sum_{\substack{i,j\ge0\\i+j\le\ell}}\frac{(\log\log y)^j}{j!}\sum_{\substack{y<q-a\le y^2\\\Omega_y(q-a)=i}}\frac{|q-a|}{\varphi(|q-a|)q}.\]
The inner sum over $q$ may again be estimated using Theorem \ref{thm:HRS} in conjunction with partial summation, yielding
\[\sum_{\substack{y<q-a\le y^2\\\Omega_y(q-a)=i}}\frac{|q-a|}{\varphi(|q-a|)q}\ll\frac{(\log\log y)^{i}}{i!\log y}.\]
Hence, the contribution to \eqref{eq:y<q-a<=y^2} from those $p\le x$ for which $\Omega_y(up+v)\le \ell$ is 
\[\ll\frac{\pi(x)}{(\log y)^2}\sum_{\substack{i,j\ge0\\i+j\le\ell}}\frac{(\log\log y)^{i+j}}{i!j!}\ll\pi(x)\cdot\frac{e^{-2Q(1/(2\log 2))\log\log y}}{\sqrt{\log\log y}}\ll\frac{\pi(x)}{(\log y)^{\eta_0}\sqrt{\log\log y}}\]
as before. This verifies \eqref{eq:y<q-a<=y^2}.

We are now able to finish the proof of Theorem \ref{thm:LMPshifted} for $3\le y\le x^{1/3}$. Observe that the number of $p\le x$ for which $up+v$ has a divisor $q-a>y$ is at most
\[\sum_{j\ge0}\left(P_{a,u,v}\big(x,y^{2^j}\big)-P_{a,u,v}\big(x,y^{2^{j+1}}\big)\right).\]
By the case $x^{1/3}<y\le x$ that we just settled, we see that the terms with $y^{2^j}>x^{1/3}$ contribute 
\[O\left(\frac{\pi(x)}{(\log y)^{\eta_0}\sqrt{\log\log y}}\right),\]
since there are at most $O(1)$ such terms. In view of \eqref{eq:y<q-a<=y^2}, we find that the contribution from the rest of the terms is
\[\ll\sum_{j\ge0}\frac{\pi(x)}{2^{\eta_0j}(\log y)^{\eta_0}\sqrt{\log\log y}}\ll\frac{\pi(x)}{(\log y)^{\eta_0}\sqrt{\log\log y}},\]
concluding the proof.
\medskip
\section*{Acknowledgments}
The author thanks Carl Pomerance for his careful reading of the manuscript and for his valuable comments. He is also grateful to Paul Pollack for helpful discussions and suggestions.

\end{document}